\documentclass{book}
\usepackage{color}
\usepackage{tmj}
\usepackage[all,ps,arc]{xy}
\usepackage{amsmath,amssymb,amsfonts}
\usepackage{graphicx}

\newcommand{\gap}{\textsf{GAP} }
\newcommand{\gapp}{\textsf{GAP}}
\newcommand{\hap}{\textsf{HAP} }
\newcommand{\happ}{\textsf{HAP}}

\newcommand{\ii}{{\bf i}}
\newcommand{\jj}{{\bf j}}

\newcommand{\Cay}{{\rm Cay}}
\newcommand{\aaa}{{\frak a}}

\begin{document}

\title{Cohomology of arithmetic groups} {Introductory computations in the cohomology of arithmetic groups}

\author{Graham~Elis} {Graham~Ellis${}^1$}

\address{${}^1$School of Mathematics, NUI Galway, Ireland
}

\email{graham.ellis@nuigalway.ie${}^1$}

\received{???} \revised{} \accepted{???}
\volumenumber{???} \volumeyear{???}

\setcounter{page}{1}

\numberwithin{equation}{section}

\maketitle

\begin{abstract}
This paper describes an approach to computer aided calculations in the 
cohomology of arithmetic groups. It  complements   existing
literature on the topic by emphasizing homotopies and perturbation techniques, rather than  
cellular subdivision, as the tools for implementing on a computer
  topological constructions that 
fail to preserve cellular structures.  
 Furthermore, it  focuses on calculating integral cohomology rather 
than just rational cohomology or cohomology at large primes.
  In particular, the paper describes and  fully implements algorithms
 for computing Hecke 
operators on the integral cuspidal cohomology of  congruence subgroups $\Gamma$ of 
  $SL_2(\mathbb Z)$, and  then partially implements versions of the 
algorithms for the special linear group   $SL_2({\cal O}_d)$ over 
various rings of imaginary quadratic integers ${\cal O}_d$. The  approach is
 also relevant for computations on  
congruence subgroups of $SL_m({\cal O}_d)$, $m\ge 2$.

\end{abstract}

\metadata{11F75}{}{cohomology of arithmetic groups, computational algebra}

\section{Introduction}

This paper aims to provide an   introduction  
 to computer aided calculations in the
cohomology of arithmetic groups up to a description, and   
partial computer implementation, of algorithms for computing Hecke
operators on the cuspidal cohomology of  congruence subgroups $\Gamma$ of
the special linear group of $2\times 2$ matrices $SL_2({\cal O})$ over
various rings of imaginary quadratic integers ${\cal O}={\cal O}_d$ as well as over the usual  integers ${\cal O}=\mathbb Z$. The approach  complements   existing
literature on the topic by emphasizing homotopies and perturbation techniques,
rather than 
cellular subdivision, as the tools for machine implementation of
 topological constructions that
fail to preserve cellular structures.
 Furthermore, we focus on calculating integral cohomology rather
than just rational cohomology or cohomology at large primes.
Implementations are available as part of the \hap package \cite{hap} for the \gap system for computational algebra \cite{GAP4}.

Section \ref{SECeichlershimura}  recalls
 some well-known motivation for studying cohomology Hecke operators.  Section \ref{SECtorsion}  recalls some  well-known motivation 
 for computing with integral, rather than rational, cohomology.
Sections \ref{SECbasic}--\ref{SECcuspidal} provide a  fully implemented 
account of how to compute  Hecke operators on the integral 
cuspidal cohomology of congruence subgroups of  $SL_2(\mathbb Z)$.   Sections \ref{SECarith}--\ref{SECch} provide a partially implemented account for congruence subgroups
of $SL_2({\cal O}_d)$ over various rings ${\cal O}_d$ of imaginary quadratic integers. The implementation is partial because the contracting homotopies of Section \ref{SECch} are not yet implemented. 
The  approach is
 also relevant for computations on
congruence subgroups of $SL_m({\cal O}_d)$, $m\ge 2$, and this is touched on in Section \ref{SECch}.

\section{The Eichler-Shimura isomorphism} \label{SECeichlershimura}

The Eichler-Shimura isomorphism 
\cite{MR89928}\cite{MR120372}
\begin{equation}
\epsilon\colon S_k(\Gamma) \oplus \overline{S_k(\Gamma)} \oplus E_k(\Gamma) 
\stackrel{\cong}{\longrightarrow} H^1(\Gamma,P_{\mathbb C}(k-2))\label{EQeichlershimura}
\end{equation}
 relates the cohomology of groups to the theory of modular forms associated to a congruence subgroup $\Gamma$ of $SL_2(\mathbb Z)$. In subsequent sections we explain how to compute with the right-hand side of the isomorphism. But first,
for completeness and for motivation,  we define the ingredients of
 the isomorphism.

 Let $N$ be a positive integer. A subgroup $\Gamma$ of $SL_2(\mathbb Z)$ is said to be a  {\em congruence subgroup} of {\em level} $N $ if it contains the kernel of the canonical homomorphism  $\pi_N\colon SL_2(\mathbb Z) \rightarrow SL_2(\mathbb Z_N)$ where $\mathbb Z_N=\mathbb Z/N\mathbb Z$.
So any congruence subgroup is of finite index in $SL_2(\mathbb Z)$.
 (To see that there exist finite index subgroups that
 are not congruence subgroups, one can use the presentation
 $SL_2(\mathbb Z) \cong \langle S, U: S^4=U^6=1, S^2=U^3\rangle$
to construct a surjective homomorphism 
$\rho\colon SL_2(\mathbb Z)\twoheadrightarrow S_7$  onto the symmetric group of degree $7$, mapping $\rho(S)= (1,2)(3,5)(4,6)$, $\rho(U)= (2,3,4)(5,6,7)$. The finite index subgroup
 $\ker \rho$  
 is not a congruence subgroup since  $S_7$ is not a quotient of $SL_2(\mathbb Z_m)$ for any $m$.)   
One congruence subgroup of particular interest is the kernel 
$\Gamma(N)=\ker \pi_N$ itself, known as the {\em principal congruence subgroup} of level $N$. A second congruence subgroup of interest is the group $\Gamma_1(N)$ 
consisting of those matrices that project to upper unitriangular matrices in
$SL_2(\mathbb Z_N)$. Another congruence subgroup of particular interest is the group  $\Gamma_0(N)$ of those matrices that project to upper triangular matrices
in $SL_2(\mathbb Z_N)$. Clearly $\Gamma(N) \le \Gamma_1(N) \le \Gamma_0(N)$.

Fix any congruence subgroup $\Gamma \le SL_2(\mathbb Z)$.
  A {\em modular form} of {\em weight} $k\ge 2$  and {\em level} 
$\Gamma $ is a complex valued  function on the upper-half plane
 $$f\colon {\frak{h}}=\{z\in \mathbb C : {\rm Re}(z)>0\} \rightarrow \mathbb C\ $$ such that  the following hold:
\begin{enumerate}
\item $\displaystyle f(\frac{az+b}{cz+d}) = (cz+d)^k f(z)$ for all $z\in \frak{h}$ and $\left(\begin{array}{ll}a&b\\ c &d \end{array}\right) \in \Gamma$, 
\item the function
	$\displaystyle (cz+d)^{-k} f(\frac{az+b}{cz+d})$ is bounded as ${\rm Im}(z) \rightarrow \infty$ for $\left(\begin{array}{ll}a&b\\ c &d \end{array}\right) \in SL_2(\mathbb Z)$,

\item $f$ is holomorphic on $\frak{h}$.

\end{enumerate}
The collection of all weight $k$ modular forms for $\Gamma$ form a vector space $M_k(\Gamma)$ over $\mathbb C$.

 A modular form $f$ is said to be a {\em cusp form} if it satisfies the following: 
\begin{enumerate}
\item[2$'$.] the function 
	$\displaystyle (cz+d)^{-k} f(\frac{az+b}{cz+d}) \rightarrow 0$ as ${\rm Im}(z)\rightarrow \infty$ for $\left(\begin{array}{ll}a&b\\ c &d \end{array}\right) \in SL_2(\mathbb Z)$. 
\end{enumerate}
The collection of all weight $k$ cusp forms for $\Gamma$ form a vector space $S_k(\Gamma)$. There is a decomposition
$$M_k(\Gamma) \cong S_k(\Gamma) \oplus E_k(\Gamma)$$
 involving a summand $E_k(\Gamma)$  known as the {\em Eisenstein
space}. 

A function  
$f\colon {\frak{h}} \rightarrow \mathbb C,z \mapsto u+{\bf i}v$ is said to be an {\em anti-holomorphic cusp form} of {\em weight} $k$ if its  complex conjugate 
 $\overline{f(z)} = u-{\bf i}v$ is a cusp form of weight $k$. The collection of all  anti-holomorphic cusp forms of weight $k$
form 
 a vector subspace $\overline{S_k(\Gamma)}$. 
See \cite{MR2289048} for further introductory details on modular forms.

On the right-hand side of (\ref{EQeichlershimura}),  the $\mathbb Z\Gamma$-module
$P_{\mathbb C}(k-2)\subset \mathbb C[x,y]$ denotes the space of homogeneous degree
$k-2$ polynomials over $\mathbb C$ with action of $\Gamma$ given by
$$\left(\begin{array}{ll}a&b\\ c &d \end{array}\right)\cdot p(x,y) = p(dx-by,-cx+ay)\ .$$
In particular $P_{\mathbb C}(0)=\mathbb C$ is the trivial module.
(In subsequent sections we  compute with the integral analogue $P_{\mathbb Z}(k-2) \subset \mathbb Z[x,y]$, to which the action of $\Gamma$ restricts.)

Each cohomology class $[c]\in H^1(\Gamma,P_{\mathbb C}(k-2))$ is represented by a function $c\colon \Gamma \rightarrow P_{\mathbb C}(k-2)$ satisfying the cocycle condition
$$ c(\gamma \gamma') = \gamma \cdot c(\gamma') + c(\gamma)$$
for all $\gamma,\gamma'\in \Gamma$. We let $Z^1(\Gamma,P_{\mathbb C}(k-2))$ denote the vector space of all such cocycles.

For any congruence subgroup $\Gamma\le SL_2(\mathbb Z)$ the 
 Eichler-Shimura map (\ref{EQeichlershimura})
is an isomorphism of vector spaces induced by the mapping
$$\epsilon\colon M_k(\Gamma) \times \overline{S_k(\Gamma)} \longrightarrow Z^1(\Gamma,P_{\mathbb C}(k-2)),\ \  (f,\overline{g}) \mapsto (c\colon \gamma \mapsto I_f({\bf i},\gamma{\bf i}) + I_{\overline g}({\bf i},\gamma {\bf i})) $$
where
$$I_f({\bf i},\gamma{\bf i}) = \int_{\bf i}^{\gamma {\bf i}} f(z)(xz +y)^{k-2}\,dz \ ,
$$
$$I_{\overline{g}}({\bf i},\gamma{\bf i}) = \int_{\bf i}^{\gamma {\bf i}} \overline{g(z)}(xz +y)^{k-2}\,d\overline{z}\ .
$$
See
\cite{wiese} for a full account of the Eichler-Shimura isomorphism.

In fact, the mapping (\ref{EQeichlershimura}) is more than an isomorphism of vector spaces. It is an isomorphism of Hecke modules: both sides admit the notion of {\em Hecke operators}, and the isomorphism preserves these operators. 
For our purposes it suffices to describe the cohomology  operators.

A congruence subgroup $\Gamma \le SL_2(\mathbb Z)$ and element 
$g\in GL_2(\mathbb Q)$  determine the finite index subgroup 
$\Gamma' = \Gamma \cap g\Gamma g^{-1} \le SL_2(\mathbb Z) $ and  homomorphisms
\begin{equation}\label{EQhomseq}
\Gamma\ \hookleftarrow\ \Gamma'\ \ \stackrel{\gamma \mapsto g^{-1}\gamma g}{\longrightarrow}\ \ g^{-1}\Gamma' g\ \hookrightarrow \Gamma\ .\end{equation}
These homomorphisms give rise to homomorphisms of cohomology groups
$$H^n(\Gamma,P_{\mathbb C}(k-2))\ \ \stackrel{tr}{\leftarrow} \ \ H^n(\Gamma',P_{\mathbb C}(k-2)) \ \ \stackrel{\alpha}{\leftarrow} \ \ H^n(g^{-1}\Gamma' g,P_{\mathbb C}(k-2)) \ \  \stackrel{\beta}{\leftarrow} H^n(\Gamma, P_{\mathbb C}(k-2))$$ 
with $\alpha$, $\beta$ functorial maps, and $tr$ the transfer map.
We define the composite 
\begin{equation}
T_g=tr \circ \alpha \circ \beta\colon H^n(\Gamma, P_{\mathbb C}(k-2)) \rightarrow H^n(\Gamma, P_{\mathbb C}(k-2)) 
\label{EQhecke}
\end{equation}
 to be the {\em Hecke operator }  determined by $g$.
The homomorphism (\ref{EQhecke}) induces homomorphisms
\begin{equation}
T_g\colon M_k(\Gamma) \rightarrow M_k(\Gamma)\ ,\ \ \ T_g\colon S_k(\Gamma) \rightarrow S_k(\Gamma)\ 
 .
\end{equation}

For $\Gamma = \Gamma_0(N)$ or $\Gamma_1(N)$ and  prime $p\ge 1$ coprime to $N$,  we  define 
the {\em Hecke operator} 
$$T_p := T_g {\rm~~~with~~~} g={\rm diag}(1,p)=\left(\begin{array}{ll} 1 &0\\ 0 &p \end{array}\right) \ .$$
This definition can be adapted to $p|N$ and composite $n$ so that it
 coincides with the classically defined Hecke operators $T_n$ for $n\ge 1$.  Further details on this description of Hecke operators can be found, for instance,
 in
\cite[Appendix by P. Gunnells]{MR2289048}.

Let $f$ be a  modular form of  weight $k\ge 2$  and  level
$\Gamma $. Suppose that the identity $f(z+1)=f(z)$ holds for all $z\in \frak{h}$. This identity certainly holds, for example, if $\Gamma=\Gamma_0(N)$ or $\Gamma=\Gamma_1(N)$. The identity can be used to establish the existence of a convergent power series
\begin{equation}
f(z) = \sum_{n=0}^\infty  a_n(e^{2\pi{\bf i}z})^n = \sum_{n=0}^\infty  a_nq^n
\label{EQpower}
\end{equation}
valid for all $z\in \frak{h}$, where $a_n$ are fixed
complex numbers and $q=e^{2\pi{\bf i}z}$.
The form $f$ is a cusp form if and only if $a_0=0$.
A non-zero
cusp form $f\in S_k(\Gamma)$ is an {\em eigenform} if it is simultaneously an eigenvector for the Hecke operators $T_n$ for all $n =1,2,3,\cdots$. An eigenform is said to be {\em normalized} if it has coefficient $a_1=1$. 

For $\Gamma=\Gamma_0(1)=SL_2(\mathbb Z)$ the vector space $S_k(SL_2(\mathbb Z))$ admits a basis of eigenforms. If $f$ is a normalized eigenform then the coefficient $a_n$ is an eigenvalue for $T_n$.
Thus, in principle,
 one can construct an  approximation to an explicit
 basis for the space $S_k(SL_2(\mathbb Z))$ of weight $k$ cusp
 forms simply by computing eigenvalues for Hecke operators.

For $\Gamma=\Gamma_0(N)$ there again exists a basis of simultaneous eigenvectors for the Hecke operators $T_n$ provided we let $n$ range over only those integers coprime to $N$.
The computation of a basis for $S_k(\Gamma_0(N))$ is a little more involved because the common eigenspaces need not be one-dimensional.
If $M$ is a positive integer dividing $N$, and if $d$ is a divisor of $N/M$,
then there is a degeneracy map $\beta_{M,d}\colon S_k(\Gamma_0(M)) \rightarrow
S_k(\Gamma_0(N))$. The {\em new subspace} of $S_k(\Gamma_0(N))$ is denoted by
$S_k^{new}(\Gamma_0(N))$ and 
 is defined to be the orthogonal complement, with respect to an inner product known as the {\em Petersson
inner product}, of the images of all maps $\beta_{ M,d }$ for $ M$ a proper divisor of $N$  and   $d | N/M$. The elements of the new subspace are called {\em newforms}. Hecke operators restrict to the space of newforms. 
 It was shown by Atkin and Lehner  \cite{MR268123} that $S_k^{new}(\Gamma_0(N))$ admits a basis of eigenforms, and that 
$$S_k(\Gamma_0(N)) = \bigoplus_{M|N}\bigoplus_{d|N/M} \beta_{M,d}(S_k^{new}(\Gamma_0(M)))\ .$$ 
 Thus, in principle,
 one can construct an  approximation to an explicit
 basis for the space $S_k(\Gamma_0(N))$ of weight $k$ cusp
 forms simply by computing eigenvalues for Hecke operators on the subspaces $S_k^{new}(\Gamma_0(M))$ of newforms. 
(An illustration is given in  Example \ref{EXeig} below for $S_2(\Gamma_0(11))$. Since this space is $1$-dimensional we have  $S_2(\Gamma_0(11))=S_2^{new}(\Gamma_0(11))$.)  
Calculating  bases of eigenforms is
 one motivation for computing Hecke operators on
 $H^1(\Gamma_0(N), P_{\mathbb C}(k-2))$.

The definition of the Hecke operators $T_g$ given in (\ref{EQhecke})  applies 
  to  
 the cohomology of any congruence subgroup $\Gamma\le SL_m(\mathbb Z)$ where $m\ge 2$, with coefficients in any finitely generated $\mathbb Z\Gamma$-module 
$\cal P$. 
 A theorem of Franke \cite{MR1603257} asserts that: (i) for suitable $\cal P$ the cohomology 
$H^\ast(\Gamma,{\cal P})$ can be directly computed in terms of certain 
automorphic forms; and (ii) there is a decomposition 
$H^\ast(\Gamma,{\cal P}) \cong H^\ast_{cusp}(\Gamma,{\cal P}) \oplus H^\ast_{eis}(\Gamma,{\cal P})$ involving a `cuspidal summand' and an `Eisenstein summand' analogous to the Eichler-Shimura isomorphism for $m=2$.
 The computation of 
 eigenvectors of Hecke operators, in this setting,   yields  
information on automorphic forms. The definition of Hecke operators applies even
more generally to congruence 
subgroups  of $SL_m({\cal O})$ with $\cal O$ the ring of integers of an algebraic number field, using elements $g\in GL_m(K)$ for the construction. 
 The {\em Bianchi case}  $m=2$ and $\cal O$ the ring of integers of an imaginary quadratic number field is considered in Sections \ref{SECarith}--\ref{SECch} but not fully implemented. 
 In this Bianchi setting  there is  a Hecke equivariant isomorphism, analogous to the Eichler-Shimura isomorphism and due to Harder \cite{MR892187}, between the space of `Bianchi modular forms' and the
first cohomology of $\Gamma$ with `appropriate coefficients'.

\section{Torsion}\label{SECtorsion}
Let $d$ be a square free integer, and let ${\cal O}_d$ denote the ring of integers of the quadratic number field $\mathbb Q(\sqrt{d})$.
Explicitly,  we have ${\cal O}_{d} = \{m+n\omega :m,n\in \mathbb Z\}$ where
$$\omega = \left\{\begin{array}{lll}
\sqrt{d}  & {\rm if~} d \equiv 2,3 {\rm ~mod~} 4\, ,\\
\frac{1+\sqrt{d}}{2}
&{\rm if~} d \equiv 1 {\rm ~mod~} 4\ .\end{array}\right.$$

Let  $\aaa \triangleleft {\cal O}_d$ be a non-zero ideal. There  is a canonical homomorphism $\pi_\aaa\colon SL_2({\cal O}_d) \rightarrow SL_2({\cal O}_d/\aaa)$.
  A
subgroup $\Gamma \le SL_2({\cal O}_d)$ is said to be a
{\em congruence subgroup} of {\em level} ${\frak a}$ if it contains $\ker \pi_\aaa$. Thus congruence subgroups are of finite index.
As above, we define 
$\Gamma({\frak a})=\ker \pi_{\frak a}$ to be the  {\em principal congruence subgroup} of level ${\frak a}$. The congruence subgroup  $\Gamma_1({\frak a})$
consists of those matrices that project to upper unitriangular matrices in
$SL_2({\cal O}_d/{\frak a})$. 
The congruence subgroup  $\Gamma_0({\frak a})$ consists of those matrices that project to upper triangular matrices
in $SL_2({\cal O}_d/{\frak a})$.

For $d>0$ the group $G=SL_2({\cal O}_{-d})$ acts on the {\em upper-half space}
$${\frak h}^3 =\{(z,t) \in \mathbb C\times \mathbb R\ |\ t > 0\} $$ 
in such a way that any point $(z,t) \in \frak{h}^3$ has finite  stabilizer group in $G$. The action is by isometries with respect to the hyperbolic metric on $\frak{h}^3$. For this metric,  geodesics  are Euclidean semi-circles of radius $0<r\le \infty$ that `meet' the complex plane $\mathbb C$ perpendicularly. 
 Any finite index subgroup $\Gamma\le G$ gives rize to a non-compact quotient orbifold $\frak{h}^3/\Gamma$ of finite volume.

The first integral homology group $H_1(\Gamma,\mathbb Z)$ is finitely generated.
Let $\Gamma^{ab}_{tors}$ be its subgroup of finite order elements. Thus  $\Gamma^{ab}_{tors}$   denotes
the maximal finite summand of $H_1(\Gamma,\mathbb Z)$. 
 Bergeron and  Venkatesh   \cite{MR3887771,MR3028790}
have conjectured relationships between
the  torsion in the integral homology of congruence subgroups $\Gamma$ and the volume of their quotient orbifold 
${\frak h}^3/\Gamma$. 
For instance, they conjecture
\begin{equation} \frac{\log |\Gamma_0(\aaa)_{tors}^{ab}|}{vol(\frak{h^3}/\Gamma_0(\aaa))} \rightarrow \frac{1}{6\pi} \label{EQbergeron} \end{equation} as the norm of the prime ideal $\aaa\triangleleft {\cal O}_{-d} $ tends to $\infty$.



Sequence (\ref{EQhomseq}), adapted to the current context, induces a composite homology homomorphism $T_g\colon H_1(\Gamma,\mathbb Z) \rightarrow H_1(\Gamma,\mathbb Z)$ associated to an element $g\in SL_2(\mathbb Q(\sqrt{-d}))$.  This 
restricts to a homomorphism $\Gamma^{ab}_{tors} \rightarrow \Gamma^{ab}_{tors}$ on the torsion part of $\Gamma^{ab}=H_1(\Gamma,\mathbb Z)$. For simplicity, let us suppose that some prime $p$ occurs with multiplicity $1$ in the prime decomposition of the order $|\Gamma^{ab}_{tors}|$. A result of P.\,Scholze \cite{MR3418533} implies, under the simplifying assumption,
 that any homology class $\alpha\in H_1(\Gamma,\mathbb Z)_p\cong \mathbb Z_p$ in the $p$-part of the homology (which is necessarily a Hecke eigenclass) 
gives rise to a representation $\rho\colon {\rm Gal}(\overline{\mathbb Q}/\mathbb Q(\sqrt{-d}))
\rightarrow GL_2(\overline{\mathbb F}_p)$ with various nice properties. Conjectures of Ash \cite{MR1788040} and others suggest a converse to Scholze's theorem.
 See 
\cite{MR3887771} for a more detailed discussion.

A computer investigation of the conjectures  of Bergeron, Venkatesh,  Ash
and others 
 is one reason for wanting  algorithms to compute with integral homology and 
cohomology.

\section{Basic computations in  $SL_2(\mathbb Z)$} \label{SECbasic}

Let $G=SL_2(\mathbb Z)$.
 The implementation of given elements of $G$,   multiplication and division of elements of $G$, and the test for equality between elements of $G$ is routine and available in all computer algebra packages. 
 The test for whether a given $2\times 2$ integer matrix lies in $G$  is an
easy test of whether certain integer equations hold, and is routine to implement. 
 The matrices
$$S=\left(\begin{array}{rr}0&-1\\ 1 &0 \end{array}\right)
{\rm~~~and~~~} T=\left(\begin{array}{rr}1&1\\ 0 &1 \end{array}\right)$$
generate $G$. It is not difficult to devise an algorithm for expressing an arbitrary integer matrix $A \in G$
 as a word
in $S$ and $T$. An
 implementation of such an algorithm underlies the functions in \hap for computing Hecke operators on the cohomology of congruence subgroups of $G$, and so we describe it in some detail. We opt for a geometric description which has the merit of being readily 
 adapted to form a key ingredient in the computation of Hecke operators for
Bianchi
 groups.  

Consider the matrix
$$U=ST = \left(\begin{array}{rr}0&-1\\ 1 &1 \end{array}\right)\ .$$ The matrices $S$ and $U$
also generate $G$. In fact we have a free presentation
 $G= \langle S,U\, |\, S^4=U^6=1, S^2=U^3 \rangle $.
The {\em cubic tree} $\cal T$ is a tree ({\em i.e.} a $1$-dimensional  contractible regular
CW-complex) with countably infinitely many edges in which each vertex has degree $3$.
We can realize  the cubic tree $\cal T$ by taking the
left cosets of ${\cal U}=\langle U\rangle$ in $G$ as vertices, and joining cosets
$x\,{\cal U} $ and $y\,{\cal U}$ by an edge if, and only if,
$x^{-1}y \in {\cal U}\, S\,{\cal U}$. Thus the
vertex $\cal U $ is joined to $S\,{\cal U} =T\,{\cal U}$,
$US\,{\cal U}= STS\,{\cal U}$ and $U^2S\,{\cal U}=T^{-1}\,{\cal U}$. The vertices of this tree
are in one-to-one
correspondence with all {\em reduced words}
 in $S$, $U$, $U^2$ 
that, apart from the identity word, end in $S$ and that don't contain the substrings $S^2$ or $U^3$.
 From this algebraic realization of the cubic tree 
we see that $G$ acts on $\cal T$  in such  a way that there is a single orbit of vertices, and a single orbit of edges; 
each vertex is stabilized by a cyclic subgroup conjugate to ${\cal U}=\langle U\rangle$ and each edge is stabilized by a cyclic subgroup conjugate to
${\cal S} =\langle S \rangle$.

\begin{figure}[h]
\centerline{\includegraphics[height=8cm]{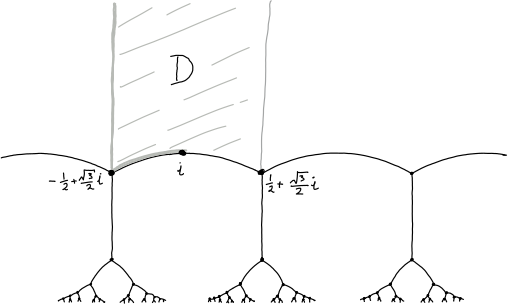}}
\caption{A portion of the cubic tree embedded in the upper-half plane $\frak h$, together with a portion of a fundamental domain $D$ for the action of $SL_2(\mathbb Z)$ on $\frak h$.}\label{FIGcubictree}
\end{figure}

Given a matrix $A \in G$ we want to describe an algorithm for producing
 a reduced word $w_A$ in $S$, $U$  and $U^2$ that represents the vertex 
$A {\cal U}$ of $\cal T$. The word $w_A$  furnishes the desired
 representation of $A$ in terms of $S$ and $T$.  For $A\in {\cal U}$ we take $w_A$ to be the empty word.  The algorithm   recursively applies 
  a procedure for determining a factorization
\begin{equation}
A=BX\label{EQfactor}
\end{equation}
of  $A\not\in {\cal U}$,  where $X\in \{S, US, U^2S\}$ and where the length of the reduced word $w_B$ is less than that of $w_A$.
One procedure for   factorization (\ref{EQfactor}) involves
 the standard action
\begin{equation}
\left(\begin{array}{ll}a &b\\ c& d\end{array}\right)\cdot z = \frac{az+b}{cz+d}
\label{EQstandardaction}
\end{equation}
of a matrix in $G$ on a point $z$ in the upper half plane $\frak{h}$. A geometric interpretation of the cubic tree  is obtained    by considering the 
singleton set and open arc
$$e^0=\{\exp({\bf i}2\pi/3)\} {\rm ~~~and~~~} 
e^1=\{\exp({\bf i}\theta) : \frac{\pi}{3} < \theta < \frac{2 \pi}{3}\}\ .$$
The union $e^0\cup e^1$ is an arc of a Euclidean unit circle,
 with one end closed
 and the other end open. The orbit of $e^0\cup e^1$ under the action of $G$ is a connected $1$-dimensional CW-complex, illustrated in Figure \ref{FIGcubictree}. The images of $e^0$ under the action are the $0$-cells, and the images of $e^1$ are the $1$-cells.
 We denote this CW-complex 
by $\cal T$  since it  
is isomorphic, as a graph, to the cubic tree 
constructed above. 
The matrix $S$ acts on $\frak h$ as inversion in the unit circle centered at $0$ followed by reflection in the imaginary axis,  $z\mapsto -1/z$.  
 The matrix $T$ acts as a  translation of one unit to the right, $z\mapsto z+1$.
The composite  $U=ST$ `rotates' through one third of a clockwise
turn  the three edges of the cubic tree touching $e^0$.

To determine  the factorization (\ref{EQfactor}) we set $z=\exp({\bf i}2\pi/3)$ 
and calculate the complex numbers $A\cdot z$ and $AX^{-1}\cdot z$ for 
$X\in \{S,US, U^2S\}$. If for one of these three choices for $X$ we find that 
the imaginary part of $AX^{-1}\cdot z$ is greater than the imaginary part of 
$A\cdot z$ then we set $B=AX^{-1}$; 
otherwise for some $X$ the absolute value of the real part of $AX^{-1}\cdot z$ is less than  the absolute value of the real part of $A\cdot z$ and the imaginary part of $AX^{-1}\cdot z$ equals that of $A\cdot z$, in which case we again set $B=AX^{-1}$.

Let $C_\ast \cal T$ denote the cellular chain complex of the CW-complex 
$\cal T$. So $C_n{\cal T}$ is the free abelian group with generators the 
$n$-cells of $\cal T$, $n=0,1$, and $d_1\colon C_1{\cal T} \rightarrow 
C_0{\cal T}$ is the boundary homomorphism.
Let ${\cal S}=\langle S\rangle$. The action of $G$ on $\cal T$ induces  
$\mathbb ZG$-module structures $C_0{\cal T}\cong \mathbb ZG\otimes_{\mathbb Z{\cal U}}\mathbb Z$ and 
$C_1{\cal T}\cong \mathbb ZG\otimes_{\mathbb Z{\cal S}}\mathbb Z^\epsilon$ where $\mathbb Z^\epsilon$ denotes the integers with non-trivial action of $\cal S$, and where 
 $d_1$ is a homomorphism of $\mathbb ZG$-modules. The
 elements of $C_0\cal T$ and $C_1\cal T$ can be expressed as $\gamma e^0$ and $\gamma e^1$
respectively, with $\gamma$ an element of the group ring $\mathbb ZG$. The boundary homomorphism is defined by $d_1(\gamma e^1) = \gamma(1-\gamma e^0)$.
The factorization (\ref{EQfactor})  can be viewed as a homomorphism 
$h_0\colon C_0{\cal T} \rightarrow C_1{\cal T}$ of free abelian groups,
recursively defined on free
generators by  
$$h_0(Ae^0) =\left\{\begin{array}{ll}
 0\, , & {\rm if~} A\in{\cal U}\, ,\\
 Be^1 + h_0(Be^0)\, ,  &{\rm if~} A\not\in {\cal U}\, .
\end{array}\right.$$ The homomorphism $h_0$ is a {\em contracting homotopy} $C_\ast {\cal T} \simeq \mathbb Z$ in the sense that
\begin{equation}
d_1h_0 = 1 - \epsilon,\ \ \ h_0d_1=1\  \label{EQcontracting}
\end{equation}
where $\epsilon\colon C_0\twoheadrightarrow H_0(C_\ast {\cal T}) \cong \mathbb Z \hookrightarrow C_0{\cal T}$ is the canonical $\mathbb Z$-linear homomorphism onto the summand $\mathbb Ze^0$ of $C_0{\cal T}$.
The homomorphism $h_0$ does not preserve the $G$-action.
The above discussion is summarized in the following.
\begin{proposition}
The $\mathbb ZG$-chain complex 
$$C_\ast {\cal T} =
 (\mathbb ZG\otimes_{\mathbb Z{\cal S}}\mathbb Z^\epsilon \stackrel{d_1}{\longrightarrow}\mathbb ZG\otimes_{\mathbb Z{\cal U}} Z)$$ and non-equivariant
contracting homotopy 
$h_\ast \colon C_\ast {\cal T}\simeq \mathbb Z$ can be implemented on a computer in
such a way that arbitrary elements $\gamma e^0\in C_0\cal T$, $\gamma e^1\in C_1\cal T$ can be expressed and their 
images $d_1(\gamma e^1)$, $h_0(\gamma e^0)$ can be uniquely determined.
\hfill{$\blacksquare$}\end{proposition}

\section{Integral cohomology of $SL_2(\mathbb Z)$}
For any group $Q$ let $R_\ast^Q$ denote some free 
$\mathbb ZQ$-resolution of the trivial module $\mathbb Z$.  In other words, 
$R_\ast^Q$ is a chain complex of free $\mathbb ZQ$-modules with 
$H_0(R^Q_\ast)\cong \mathbb Z$, $H_n(R^Q_\ast)=0$ for $n>0$.
The cohomology of $Q$ with coefficients in the trivial $Q$-module $\mathbb Z$ is defined as
$$H^n(Q,\mathbb Z) = H^n({\rm Hom}_{\mathbb ZQ}(R^Q_\ast,\mathbb Z))\ .$$
A free resolution $R^Q_\ast$ always admits a contracting homotopy $h_\ast\colon
 R^Q_\ast \simeq \mathbb Z$  consisting of
 a sequence of $\mathbb Z$-linear homomorphisms 
$h_n\colon R_n^Q\rightarrow R_{n+1}^Q$ for $n\ge 0$
 satisfying $d_{n+1}h_{n} +h_{n-1}d_n =1\  (n>0)$, $d_1h^0=1-\epsilon$
 where $\epsilon\colon R^Q_0\twoheadrightarrow H_0(R^Q_\ast) \cong \mathbb Z \hookrightarrow R_0^Q$ is the canonical $\mathbb Z$-linear homomorphism onto the summand $\mathbb Z$ of $R^Q_0$.

Many theoretical
constructions in the cohomology of groups  
 involve  repeated use of the following  element of choice.
\begin{quote} {\bf Element of choice:} Given $x\in \ker (d_n\colon R^Q_n\rightarrow R^Q_{n-1})$  choose an element $\tilde x\in R^Q_{n+1}$ such that $d_{n+1}(\tilde x) = x$.
\end{quote}
 If an algorithmic
 formula for a contracting homotopy on $R^Q_\ast$ is to hand then the choice can be made algorithmic: one simply
chooses $\tilde x = h_n(x)$.
 
For a cyclic group $Q =\langle x : x^q=1\rangle$  one can choose $R^Q_n=\mathbb ZQ$ for $n\ge 0$, and $d_{2n-1}(1)=(x-1)$, $d_{2n}(1)=(1+x+x^2+\cdots+x^{q-1})$
for $n>0$.
A contracting homotopy is given by $h_{2n}(x^k)=1+x+x^2+\cdots+x^{k-1}$,
$h_{2n+1}(x^{q-1})=1$, $h_{2n+1}(x^\ell)=0$ for $\ell\ne q-1$.

For the  group $G=SL_2(\mathbb Z)$, and specific cyclic subgroups
 ${\cal U}=\langle U\rangle$,
${\cal S}=\langle S\rangle$ we have 
$$H_0(R^{\cal U}_\ast \otimes_{\mathbb Z{\cal U}} \mathbb ZG) \cong \mathbb ZG\otimes_{\mathbb Z{\cal U}} \mathbb Z =C_0{\cal T} {\rm ~~~and~~~}
H_0((R^{\cal S}_\ast \otimes_{\mathbb Z } \mathbb Z^\epsilon) \otimes_{\mathbb Z{\cal S} } \mathbb ZG) \cong \mathbb ZG\otimes_{\mathbb Z{\cal S}} \mathbb Z^\epsilon =C_1{\cal T}\ .$$
The boundary homomorphism $d_1\colon C_1{\cal T}\rightarrow C_0{\cal T}$ thus induces a chain homomorphism 
\begin{equation}
d_\ast^h\colon (R^{\cal S}_\ast \otimes_{\mathbb Z}\mathbb Z^\epsilon) \otimes_{\mathbb Z{\cal S}} \mathbb ZG \longrightarrow R^{\cal U}_\ast \otimes_{\mathbb Z{\cal U}} \mathbb ZG
\label{EQdouble}\end{equation}
between free $\mathbb ZG$-chain complexes. The superscript on $d_\ast^h$ stands for `horizontal'. We  regard (\ref{EQdouble}) as a double complex, and let $R^G_\ast$ denote its total complex. Explicitly $R^G_n=D_{1,n-1} \oplus D_{0,n}$ where
$$D_{0,n} = (R^{\cal U}_{n} \otimes_{\mathbb Z{\cal U}} \mathbb ZG)\cong \mathbb ZG\ , \ \ \
D_{1,n-1} = ((R^{\cal S}_{n-1} \otimes_{\mathbb Z} \mathbb Z^\epsilon) \otimes_{\mathbb Z{\cal S}} \mathbb ZG) \cong \mathbb ZG
$$
with $R^{\cal S}_{-1}=0$. The boundary homomorphism on $R^G_\ast$ is
$$d_n\colon D_{1,n-1}\oplus D_{0,n} \rightarrow D_{1,n-2}\oplus D_{0,n-1}, x\oplus y\mapsto d^v_{n-1}(x) + (-1)^nd_{n-1}^h(x) +d^v_n(y) $$
where the `vertical' homomorphisms $d^v_n$ are induced by the boundary maps on $R^{\cal U}_\ast$ and $R^{\cal S}_\ast$.
 The spectral sequence of a double complex together with the exactness of
 the complexes $(D_{1,\ast},d^v_\ast)$ and $(D_{0,\ast},d^v_\ast)$ imply 
that the free $\mathbb ZG$-chain complex $R^G_\ast$ is a resolution of $\mathbb Z$. A contracting homotopy  $h_\ast\colon R^G_\ast \rightarrow R^G_{\ast+1}$ can be constructed from  contracting homotopies
$h_\ast^v\colon R^{\cal U}_\ast\rightarrow R^{\cal U}_{\ast+1}$, $h_\ast^v\colon R^{\cal S}_\ast\rightarrow R^{\cal S}_{\ast+1}$,
$h_0^h\colon C_0{\cal T}\rightarrow C_1{\cal T}$ using the formula
$$h_n(x\oplus y) = h_{n-1}^v(x) \oplus \{(-1)^nh^v_nd^h_nh^v_{n-1}(x) + h^v_n(y)\}\ ,$$
$$h_0(y) = h_0^h(y)\oplus \{h^v_0(y)-h^v_0d^h_0 h^h_0(y)\}\ .$$
In these formulas $h^v_n,h^h_0$ denote the maps induced by tensoring.
In summary,  we have established the following.
\begin{proposition}\label{PROPtwo}
Let $G=SL_2(\mathbb Z)$. A  free $\mathbb ZG$-resolution $R^G_\ast$ of $\mathbb Z$ and contracting homotopy $h_\ast \colon R^G_\ast \simeq \mathbb Z$ can be implemented on a computer, with $R^G_0=\mathbb ZG$, $R^G_n=\mathbb ZG\oplus \mathbb ZG$ for $n\ge 1$. 
Arbitrary elements $w \in R^G_n$  can be expressed and their
images $d_n(w)$, $h_n(w)$ can be uniquely determined.
\hfill{$\blacksquare$}\end{proposition}

\section{Integral cohomology of congruence subgroups of $SL_2(\mathbb Z)$}
Let $G=SL_2(\mathbb Z)$ and let $\Gamma$ denote a congruence subgroup for which we can algorithmically test membership $A\in \Gamma$ for any matrix $A$ in $G$.
 For instance, $\Gamma$ could be one of
 the congruence subgroups $\Gamma(N), \Gamma_1(N), \Gamma_0(N)$ of level $N$.

Let $K=\Cay(G)$  be the Cayley graph of $G$ with respect to the generators $S,U$. The vertices of $K$ are the elements of $G$ and there is a single edge between vertices $A,A'\in G$ if, and only if, $A^{-1}A'\in \{S,U\}$ or
$A'^{-1}A\in \{S,U\}$.
We can choose some vertex $v_0$ in $K$ and, using the membership test for $\Gamma$,
 perform a breadth first search of the graph $K$ in order
to construct some connected subgraph $D$ of $K$ that contains $v_0$ and 
that is maximal with respect to the property that the vertices of $D$ belong to distinct orbits under the action of $\Gamma$. 
An edge of $K$ with precisely one boundary vertex in $D$  corresponds to an element of $\Gamma$, and the
 collection of such edges determines a finite generating set for $\Gamma$.
This generating set likely contains many redundant generators, and we can try to form a smaller generating set by searching for obvious   
  redundancies.
\begin{figure}[h]
\centerline{\includegraphics[height=4cm]{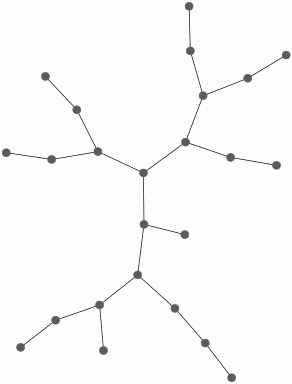} \ \ \ \ \includegraphics[height=6cm]{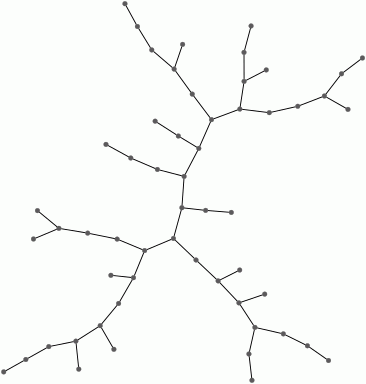} }
\caption{Connected subtrees $D$ for the congruence subgroups
$\Gamma(6)$ (left) and $\Gamma_0(39)$ (right). }\label{FIGmaxtree}
\end{figure}

Figure \ref{FIGmaxtree} (right) shows a maximal subtree $D$ for the action of 
$\Gamma_0(39)$ on the Cayley graph $K=\Cay(G)$. It has $56$ vertices, indicating that $\Gamma_0(39)$ is of index $56$ in $G$. The subtree  yields a generating set for $\Gamma_0(39)$ which, after elimination of obvious redundancies, consists of $18$ generators. 
The vertices of  $D$ represent a 
 transversal of $\Gamma$ in $G$ consisting of $|G:\Gamma|$ coset representatives. In examples such as this, where  $D$ is fairly small, it is practical to determine
 the transversal element $t$ representing an arbitrary element $g\in G$ by naively iterating over the transversal until the transversal element $t$  satisfying $tg^{-1}\in \Gamma$ is found. This provides a permutation action $G\rightarrow S_{|G:\Gamma|}$ of $G$ on the transversal.

For larger index $|G:\Gamma|$ it can be more efficient  
 to work with a connected graph  on which $G$ acts so that
 the vertices have non-trivial stabilizer groups.   
In particular, for $\Gamma=\Gamma(N)$ with $N\ge 3$ we  can take $K=\cal T$ to be 
the cubic tree so that each vertex has  stabilizer group in $G$ of  order $6$. 
The action of $\Gamma$  on $\cal T$ factors through an action of $\overline {\Gamma}=\Gamma(N)/\langle S^2\rangle \le PSL_2(\mathbb Z)$. The group $\overline {\Gamma}$ acts  freely on  $\cal T$.  
 We can thus use the above method to find a generating set for $\overline{\Gamma}$ and lift it to a generating set for $\Gamma$.   
Figure \ref{FIGmaxtree} (left) shows a  maximal subtree $D$  for the action of the principal congruence subgroup $\Gamma(6)$  on the cubic tree $\cal T$; it has $24$ vertices, indicating that $\Gamma(6)$ is of index $144=6\times 24$ in $G$; the subtree yields a generating set of $13$ generators for $\Gamma(6)$.

Any
 free $\mathbb ZG$-resolution $R^G_\ast$ is also a free
$\mathbb Z\Gamma$-resolution, where ${\rm rank}_{\mathbb Z\Gamma}(R^G_n) = |G:\Gamma| \times {\rm rank}_{\mathbb ZG}(R^G_n)$.
 In light of Proposition \ref{PROPtwo}, and the use of  contracting homotopies to make algorithmic the element of choice in constructing Hecke operators, we have established the following.

\begin{proposition}
Let $\Gamma$ be  a congruence subgroup of $G$ with an algorithmic membership test and let $g\in G$. We can  
implement the Hecke operator $T_g\colon H^m(\Gamma, \mathbb Z) \rightarrow H^m(\Gamma, \mathbb Z)$ on a computer.
\hfill{$\blacksquare$}\end{proposition} 

\begin{example}
	The following \hap commands compute the  Hecke operators $T_{g_n}\colon H^1(\Gamma(6),\mathbb Z) \rightarrow H^1(\Gamma(6),\mathbb Z)$ on weight $k=2$ forms for $g_n={\rm diag}(1,n), n=2,5$, and confirm that $T_{g_2}T_{g_5}=T_{g_5}T_{g_2}$.
\begin{verbatim}
gap> gamma:=HAP_PrincipalCongruenceSubgroup(6);;
gap> n:=2;;Tg2:=HeckeComponentWeight2(gamma,n,1);;
gap> M2:=HomomorphismAsMatrix(Tg2);;Display(M2);
[ [   1,   0,   0,   0,   0,   0,   0,   0,   0,   0,   0,   1,   0 ],
  [   0,   0,   0,   0,   2,   0,   1,   0,   0,   0,  -1,   0,   0 ],
  [   0,   0,   0,   2,   0,   0,   0,   0,   0,   0,   0,   0,  -1 ],
  [   0,   0,   0,   0,   2,   0,   1,   0,   0,   0,  -1,   0,   0 ],
  [   0,   0,   0,   2,   0,   0,   0,   0,   0,   0,   0,   0,  -1 ],
  [   0,   1,   0,   0,   0,   0,   1,   1,   0,   0,  -1,   0,   0 ],
  [   0,   0,   0,   0,   0,   0,   1,   0,   0,   0,   1,   0,  -1 ],
  [   0,   0,   1,   0,   0,   1,   0,   0,   0,   0,   0,   0,  -1 ],
  [   0,   0,   0,   0,   0,   0,   0,   0,   1,   1,   0,   0,  -1 ],
  [   0,   0,   0,   0,   0,   0,   0,   0,   0,   0,   0,   0,   0 ],
  [   0,   0,   0,   0,   0,   0,   1,   0,   0,   0,   1,   0,  -1 ],
  [   0,   0,   0,   0,   0,   0,   1,   0,   0,   0,   1,   0,  -1 ],
  [   0,   0,   0,   0,   0,   0,   0,   0,   0,   0,   0,   0,   0 ] ]
gap> n:=5;;Tg5:=HeckeComponentWeight2(gamma,n,1);;
gap> M5:=HomomorphismAsMatrix(Tg5);;Display(M5);
[ [   6,   0,   0,   0,   0,   0,   1,   0,   0,   0,  -1,   0,  -1 ],
  [   0,   0,   6,   0,   0,   0,  -1,   0,   0,   0,   1,   0,  -2 ],
  [   0,   6,   0,   0,   0,   0,   2,   0,   0,   0,  -2,   0,   1 ],
  [   0,   0,   0,   0,   6,   0,   3,   0,   0,   0,  -3,   0,   0 ],
  [   0,   0,   0,   6,   0,   0,   0,   0,   0,   0,   0,   0,  -3 ],
  [   0,   0,   0,   0,   0,   0,   4,   6,   0,   0,  -4,   0,  -1 ],
  [   0,   0,   0,   0,   0,   0,   3,   0,   0,   0,   3,   0,  -3 ],
  [   0,   0,   0,   0,   0,   6,   1,   0,   0,   0,  -1,   0,  -4 ],
  [   0,   0,   0,   0,   0,   0,  -2,   0,   6,   0,   2,   0,  -4 ],
  [   0,   0,   0,   0,   0,   0,   2,   0,   0,   6,  -2,   0,  -2 ],
  [   0,   0,   0,   0,   0,   0,   3,   0,   0,   0,   3,   0,  -3 ],
  [   0,   0,   0,   0,   0,   0,  -1,   0,   0,   0,   1,   6,   1 ],
  [   0,   0,   0,   0,   0,   0,   0,   0,   0,   0,   0,   0,   0 ] ]
gap> M2*M5=M5*M2;
true
\end{verbatim}
\end{example}

For $\Gamma=\Gamma_0(N)$ a permutation action $G\rightarrow S_{|G:\Gamma_0(N)|}$ can be constructed more efficiently using the following well-known result which we recall from \cite{MR2467560}, particularly  when $N$ is prime.
See  \cite[Proposition 2.2.2]{MR1628193} for a proof.

Let
$$\mathbb P^1(\mathbb Z_N) =\{(a:b) \,|\,a,b\in\mathbb Z_N, gcd(a,b,N)=1)\}/\sim$$
where $(a:b)\sim(a':b')$ if there is a unit $u$ in $\mathbb Z_N$
such that $a=ua', b=ub'$.
\begin{proposition}\label{PROPcre}
There is an equivariant bijection between $\mathbb P^1(\mathbb Z_N)$ and the right cosets of $\Gamma_0(N)$ in $SL_2(\mathbb Z)$, which sends a coset representative
$\left(\begin{array}{ll}a &b\\c&d\end{array}\right)$ to the class of $(c:d)$ in $\mathbb P^1(\mathbb Z_N)$.
\end{proposition}
\section{Simple homotopy collapses}
An obvious bottleneck in the above approach to cohomology calculations
 for 
congruence subgroups $\Gamma\le G=SL_2(\mathbb Z)$ is the rank of the 
modules in the $\mathbb ZG$-resolution $R^G_\ast$ when considered as free $\mathbb Z\Gamma$-modules. 

\begin{example} \label{EXshc}
The congruence subgroup $\Gamma=\Gamma_0(1000)$ is of index $1800$ in $G$. Let 
$S^\Gamma_\ast$ be the free $\mathbb Z\Gamma$-resolution of $\mathbb Z$ obtained from the resolution $R^G_\ast$ of Proposition \ref{PROPtwo}
 by considering each $R^G_n$ as a $\mathbb Z\Gamma$-module.
Then ${\rm rank}_{\mathbb Z\Gamma}\,S^\Gamma_0 =1800$  
and ${\rm rank}_{\mathbb Z\Gamma}\,S^\Gamma_n =3600$ for $n\ge 1$.
To calculate, for example, the homology group $H_5(\Gamma,\mathbb Z)=H_5(C_\ast)$ directly from the chain complex
$C_\ast=S^\Gamma_\ast\otimes_{\mathbb Z\Gamma}\mathbb Z$
would involve an application of the Smith Normal Form algorithm to a boundary matrix of dimensions $3600\times 3600$, and such an application would challenge the efficient implementation of the SNF algorithm available
in \gap. To avoid this challenge we could try to find a chain
homotopy equivalence $S^\Gamma_\ast \simeq T^\Gamma_\ast$ 
with $T^\Gamma_\ast$  a smaller chain complex of free $\mathbb Z\Gamma$-modules and 
compute the required homology from the chain complex
$D_\ast=T^\Gamma_\ast\otimes_{\mathbb ZG}\mathbb Z$; alternatively we could try to compute a chain homotopy equivalence $C_\ast\simeq D_\ast$ directly.
The following \hap commands use the latter approach to compute 
$H_5(\Gamma,\mathbb Z) =\mathbb Z_5$ in a way that involves an application of the SNF algorithm to a matrix of dimensions $302\times 302$.
\begin{verbatim}
gap> gamma:=HAP_CongruenceSubgroupGamma0(1000);;
gap> R:=ResolutionSL2Z(1,6);;
gap> S:=ResolutionFiniteSubgroup(R,gamma);;
gap> C:=TensorWithIntegers(S);;
gap> List([0..5],C!.dimension);
[ 1800, 3600, 3600, 3600, 3600, 3600 ]
gap> D:=ContractedComplex(C);;
gap> List([0..5],D!.dimension);
[ 1, 302, 302, 302, 302, 302 ]
gap> Homology(D,5);
[ 2 ]
\end{verbatim}
\end{example}

To explain how the homotopy equivalence $C_\ast \simeq D_\ast$ was 
constructed in Example \ref{EXshc} let us consider an arbitrary chain complex 
$C_\ast$ of free $\Lambda$-modules $C_n$ where $\Lambda$ is an associative 
(but not necessarily commutative) ring with identity. The examples we have in 
mind are $\Lambda=\mathbb Z$ and $\Lambda=\mathbb Z\Gamma$.
Let us denote the free generators of $C_n$ by $e^n_1, \cdots, e^n_k$, $k={\rm rank}_{\Lambda}\,C_n$. The boundary homomorphism is given by $d_n(e^n_i) = \lambda_{i1}e^{n-1}_1 + \cdots +\lambda_{i\ell}e^{n-1}_\ell$ with $\ell = {\rm rank}_{\Lambda}\, C_{n-1}$.
 Suppose that for some particular generator $e^n_i$ one of the coefficients 
$\lambda_{ij}$ is a unit in $\Lambda$. Let $\langle e^n_i, d_n(e^n_i)\rangle$ denote the sub $\Lambda$-chain complex generated by $ e^n_i$ and $ d_n(e^n_i)$. 
Since one of the coefficients is a unit, this sub chain complex has trivial homology, 
and the quotient chain complex
 $C'_\ast =C_\ast/\langle e^n_i, d_n(e^n_i)\rangle$ is a chain complex of free
$\Lambda$-modules. It follows from the exact homology sequence of a short exact sequence of chain complexes that the
  quotient chain map 
$C_\ast \twoheadrightarrow C'_\ast$ is a quasi-isomorphism and thus  homotopy equivalence of chain complexes.
 We say that $C'_\ast$ is obtained from $C_\ast$ by a {\em simple homotopy collapse} and write $C_\ast \searrow C'_\ast$. We can search, recursively, for a sequence of simple homotopy collapses $C_\ast \searrow C'_\ast \searrow C''_\ast \searrow C'''_\ast \searrow \cdots \searrow D_\ast$
and use $D_\ast$ in place of $C_\ast$ in cohomology computations.
 Example \ref{EXshc} illustrates this technique for $\Lambda=\mathbb Z$. The next example illustrates the technique for $\Lambda=\mathbb Z\Gamma$.

\begin{example}
The congruence subgroup $\Gamma=\Gamma_0(50)$ is of index $90$ in $G=SL_2(\mathbb Z)$. Let
$S^\Gamma_\ast$ be the free $\mathbb Z\Gamma$-resolution  obtained from $R^G_\ast$ by restricting the action.
The following \hap commands construct a chain homotopy equivalence $S^\Gamma_\ast \simeq T^\Gamma_\ast$ with ${\rm rank}_{\mathbb Z\Gamma}\,T^\Gamma_0=1$, ${\rm rank}_{\mathbb Z\Gamma}\,T^\Gamma_n=17$ for $n\ge 1$ and use $T^\Gamma_\ast$ to compute 
$H^1(\Gamma,P_{\mathbb Z}(4))=\mathbb Z_{2}\oplus \mathbb Z_4\oplus \mathbb Z_{120} \oplus \mathbb Z^{174}$ and
$H^5(\Gamma,P_{\mathbb Z}(4))=\mathbb Z_{2}^{77}$.
\begin{verbatim}
gap> gamma:=HAP_CongruenceSubgroupGamma0(50);;
gap> R:=ResolutionSL2Z(1,6);;
gap> S:=ResolutionFiniteSubgroup(R,gamma);;
gap> List([0..5],S!.dimension);
[ 90, 180, 180, 180, 180, 180 ]
gap> T:=ContractedComplex(S);;
gap> List([0..5],T!.dimension);
[ 1, 17, 17, 17, 17, 17 ]
gap> P:=HomogeneousPolynomials(gamma,4);;
gap> C:=HomToIntegralModule(T,P);;
gap> Cohomology(C,1);
[ 2, 4, 120, 0, 0, 0, 0, 0, 0, 0, 0, 0, 0, 0, 0, 0, 0, 0, 0, 0, 0, 0, 0, 0, 
  0, 0, 0, 0, 0, 0, 0, 0, 0, 0, 0, 0, 0, 0, 0, 0, 0, 0, 0, 0, 0, 0, 0, 0, 0, 
  0, 0, 0, 0, 0, 0, 0, 0, 0, 0, 0, 0, 0, 0, 0, 0, 0, 0, 0, 0, 0, 0, 0, 0, 0, 
  0, 0, 0 ]
gap> Cohomology(C,5);
[ 2, 2, 2, 2, 2, 2, 2, 2, 2, 2, 2, 2, 2, 2, 2, 2, 2, 2, 2, 2, 2, 2, 2, 2, 2, 2, 
  2, 2, 2, 2, 2, 2, 2, 2, 2, 2, 2, 2, 2, 2, 2, 2, 2, 2, 2, 2, 2, 2, 2, 2, 2, 2, 
  2, 2, 2, 2, 2, 2, 2, 2, 2, 2, 2, 2, 2, 2, 2, 2, 2, 2, 2, 2, 2, 2, 2, 2, 2 ]
\end{verbatim}
\end{example}

\section{Cuspidal cohomology of congruence subgroups of $SL_2(\mathbb Z)$} \label{SECcuspidal}
The action (\ref{EQstandardaction}) of $G=SL_2(\mathbb Z)$ on the upper-half plane $\frak h$ has a fundamental domain
$$\begin{array}{lclcl}
D&=&\{z\in {\frak h} : |z|>1, |{\rm Re}(z)|<\frac{1}{2}\}
&\cup &\{z\in {\frak h} : |z|\ge 1, {\rm Re}(z)=-\frac{1}{2}\}\medskip\\
&&&\cup &\{z\in {\frak h} : |z|=1, -\frac{1}{2} \le {\rm Re}(z) \le 0\} 
\end{array}
$$
shown in Figure \ref{FIGcubictree}.
The action factors through an action of
$PSL_2(\mathbb Z) =SL_2(\mathbb Z)/\langle
\left(\begin{array}{rr}-1&0\\ 0 &-1 \end{array}\right)\rangle$.  The images of $D$ under the action of $PSL_2(\mathbb Z)$
cover the upper-half plane, and any two images have at most a single point in common. The possible common points are in the orbit of the bottom left-hand corner point 
$-\frac{1}{2}+{\bf i}\frac{\sqrt{3}}{2}$ which is stabilized by 
$\cal U$, or in the orbit of  the bottom middle point ${\bf i} $
which is stabilized by 
$\cal S$.

A congruence subgroup $\Gamma$
has a `fundamental domain' $D_\Gamma$ equal to a union of finitely many
copies of $D$, one copy for each coset in $\Gamma\setminus SL_2(\mathbb Z)$.
The
quotient space $X=\Gamma\setminus {\frak h}$ is not compact, and can be
compactified in several ways. We are interested in the Borel-Serre
compactification.
This is a space $X^{BS}$ for which there is an inclusion
$X\hookrightarrow X^{BS}$ that is a homotopy equivalence.
One defines the  {\em boundary} $\partial X^{BS} = X^{BS} - X$
 and uses the inclusion $\partial X^{BS} \hookrightarrow X^{BS} \simeq X$
 to define the cuspidal cohomology group, over the ground ring $\mathbb C$, as
$$H_{cusp}^n(\Gamma,P_{\mathbb C}(k-2)) = \ker (\ H^n(X,P_{\mathbb C}(k-2)) \rightarrow
H^n(\partial X^{BS},P_{\mathbb C}(k-2)) \ ).$$
Strictly speaking, this is the definition of  {\em interior cohomology}
$H_!^n(\Gamma,P_{\mathbb C}(k-2))$ which in general contains the
 cuspidal cohomology as a subgroup. However, for congruence subgroups of
$SL_2(\mathbb Z)$ there is equality
$H_!^n(\Gamma,P_{\mathbb C}(k-2)) = H_{cusp}^n(\Gamma,P_{\mathbb C}(k-2))$.

Working over $\mathbb C$ has the advantage of
avoiding the technical issue that $\Gamma $ does not necessarily act freely on 
${\frak h}$ since there may be points with finite
cyclic stabilizer groups in $SL_2(\mathbb Z)$. But it has the disadvantage of 
losing  information about torsion in cohomology. We address the issue
by working with a contractible CW-complex
$\tilde X^{BS}$ on which $\Gamma$ acts freely, and $\Gamma$-equivariant inclusion
$\partial \tilde X^{BS} \hookrightarrow \tilde X^{BS}$. The definition of cuspidal cohomology that we use, which coincides with the above definition when working over $\mathbb C$, is
\begin{equation}
H_{cusp}^n(\Gamma,A) = \ker (\ H^n({\rm Hom}_{\, \mathbb  Z\Gamma}(C_\ast(\tilde X^{BS}), A)\,   ) \rightarrow
H^n(\ {\rm Hom}_{\, \mathbb  Z\Gamma}(C_\ast(\partial \tilde X^{BS}), A)\, \ ).
\label{EQcusp}
\end{equation}

\begin{figure}[h]
\centerline{\includegraphics[height=5cm]{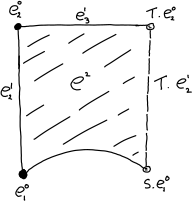}}
\caption{Fundamental domain for the action of $PSL_2(\mathbb Z)$ on $X^{BS}$.}\label{FIGfd}
\end{figure}

The compact CW-complex $X^{BS}$ 
is described by the CW-structure on the fundamental domain for its action of $G$ shown in Figure \ref{FIGfd}  
and the cell stabilizer groups
$Stab(e^0_1)={\cal U} \cong C_6$, 
$Stab(e^0_2)=\langle U^3\rangle \cong C_2$,
$Stab(e^1_1)={\cal S} \cong C_4$,
$Stab(e^1_2)=\langle U^3\rangle \cong C_2$,
$Stab(e^1_3)=\langle U^3\rangle \cong C_2$,
$Stab(e^2)=\langle U^3\rangle \cong C_2$.
The cellular chain comlex $C_\ast X^{BS}$ is a complex of $\mathbb ZG$-modules of the form
$$0\longrightarrow \mathbb ZG\otimes_{C_2} \mathbb Z \longrightarrow
\mathbb ZG\otimes_{C_4} \mathbb Z^\epsilon \,\oplus\,
\mathbb ZG\otimes_{C_2} \mathbb Z \,\oplus\,
\mathbb ZG\otimes_{C_2} \mathbb Z
\longrightarrow
\mathbb ZG\otimes_{C_6} \mathbb Z
\,\oplus\, \mathbb ZG\otimes_{C_2} \mathbb Z\ .$$
The process of using resolutions for cell stabilizer groups to convert the contractible $\mathbb ZG$-complex $C_\ast {\cal T}$ into a free $\mathbb ZG$-resolution $R^G_\ast$ can be adapted to the current setting.
Resolutions for cell stabilizers can be combined with the $\mathbb ZG$-complex 
$C_\ast X^{BS}$ to produce a free $\mathbb ZG$-resolution $C_\ast\tilde X^{BS}$.
The construction uses a perturbation technique of C.T.C.\,Wall \cite{MR0178046} 
and explicit formulas for the construction in terms of contracting homotopies
can be found in \cite{MR2270569}. The following is a summary of the construction.
\begin{proposition}\cite{MR2270569} \label{PROPwall}
Let $X$ be any contractible CW-complex on which some group $G$ acts in a way that permutes cells. 
Suppose that for $n\ge 0$ there are finitely many orbits of $n$-cells 
represented by $e^n_1, e^n_2, \cdots e^n_{k_n}$.
 Let $G^{e^n_i}$ denote the subgroup of $G$  stabilizing $e^n_i$. 
Suppose that we have  free $\mathbb ZG^{e^n_i}$-resolutions 
$R^{G^{e^n_i}}_\ast$ of $\mathbb Z$. Then there is a free $\mathbb ZG$-resolution $R^G_\ast$ of 
$\mathbb Z$ with $$R^G_n = \bigoplus_{p+q=n, p,q\ge 0} (R^{G^{e^p_i}}_q\otimes_{\mathbb Z} \mathbb Z^{\epsilon^p_i}) \otimes_{\mathbb ZG^{e^p_i}}\mathbb ZG $$
where $\mathbb Z^{\epsilon^p_i}$ denotes the integers with some action of $G^{e^p_i}$.
 An explicit formula for the boundary homomorphism $d_n\colon R^G_n\rightarrow R^G_{n-1}$ can be given in terms of the boundary homomorphism on $C_\ast X$ and the boundary homomorphisms and contracting homotopies for the resolutions of the stabilizer groups. An explicit formula for a contracting homotopy $h_n\colon R^G_n\rightarrow R
^G_{n+1}$ can also be given if, in addition, we have an explicit formula for a contracting homotopy $h_\ast\colon C_\ast X\simeq \mathbb Z$.
\hfill{$\blacksquare$}\end{proposition}
The free resolution $C_\ast\tilde X^{BS}$ is of the form
${\rm rank}_{\mathbb ZG}(C_0\tilde X^{BS})=2$, 
${\rm rank}_{\mathbb ZG}(C_1\tilde X^{BS})=5$, and 
${\rm rank}_{\mathbb ZG}(C_n\tilde X^{BS})=6$ for $n\ge 2$.
Having constructed $C_\ast \tilde X^{BS}$ it is routine to implement the
definition \ref{EQcusp} of cuspidal cohomology.

\begin{example}
The following \hap commands compute $H^1_{cusp}(\Gamma_0(39),P_{\mathbb Z}(2))\cong \mathbb Z^{24}$.
\begin{verbatim}
gap> gamma:=HAP_CongruenceSubgroupGamma0(39);;
gap> k:=4;; deg:=1;; c:=CuspidalCohomologyHomomorphism(gamma,deg,k);;
gap> AbelianInvariants(Kernel(c));
[ 0, 0, 0, 0, 0, 0, 0, 0, 0, 0, 0, 0, 0, 0, 0, 0, 0, 0, 0, 0, 0, 0, 0, 0 ]
\end{verbatim}
\end{example}

\begin{example}\label{EXeig}
The following \hap commands establish that $S_2(\Gamma_0(11))\cong \mathbb C$ is 1-dimensional with basis eigenform
$$f = q -2q^2 -q^3 +2q^4 +q^5 +2q^6 -2q^7  -2q^9 -2q^{10} + \cdots\ .$$
\begin{verbatim}
gap> gamma:=HAP_CongruenceSubgroupGamma0(11);;
gap> AbelianInvariants(Kernel(CuspidalCohomologyHomomorphism(gamma,1,2)));
[ 0, 0 ]
gap> gamma:=HAP_CongruenceSubgroupGamma0(11);;
gap> for n in [1,2,3,5,7] do
> Display(HomomorphismAsMatrix(HeckeComponentWeight2(gamma,n,1)));;
> od;
[ [  1,  0,  0 ],   [ [   3,  -4,   4 ],   [ [   4,  -4,   4 ],
  [  0,  1,  0 ],     [   0,  -2,   0 ],     [   0,  -1,   0 ],
  [  0,  0,  1 ] ]    [   0,   0,  -2 ] ]    [   0,   0,  -1 ] ]
\end{verbatim}
\begin{verbatim}
[ [   6,  -4,   4 ],    [ [   8,  -8,   8 ],
  [   0,   1,   0 ],      [   0,  -2,   0 ],
  [   0,   0,   1 ] ]     [   0,   0,  -2 ] ]
\end{verbatim}
\end{example}
As explained in \cite{MR2467560}, for a normalized eigenform $f=1 + \sum_{s=2}^\infty a_sq^s \in S_k(\Gamma_0(N))$
the coefficients $a_s$ with $s$ a composite integer can be expressed in terms of the coefficients $a_p$ for prime $p$.
If $r,s$ are coprime then $a_{rs} =a_ra_s$.
If $p$ is a prime that is not a divisor of the level $N$ of $\Gamma$ then
$a_{p^m} =a_{p^{m-1}}a_p - p a_{p^{m-2}}.$
If the prime $ p$ divides $N$ then $a_{p^m} = (a_p)^m$. It thus suffices to compute the coefficients $a_p$ for prime integers $p$ only.

See Stein's paper \cite{MR2467560} for other techniques for
 computing  Fourier expansions of 
 classical modular forms,  in particular techniques using Manin symbols.

\section{Integral cohomology of $SL_2(\mathbb Z[\ii])$}\label{SECarith}

The group $SL_2(\mathbb C)$ acts on  
the  {\em upper-half space}
$${\frak h}^3 =\{(z,t) \in \mathbb C\times \mathbb R\ |\ t > 0\}  $$
in a  well-known fashion.
To describe the action we introduce the symbol $\jj$ satisfying $\jj^2=-1$, 
$\ii\jj=-\jj\ii$ where $\ii=\sqrt{-1}$, and write $z+t\jj$ instead of $(z,t)$. The action is given by
$$\left(\begin{array}{ll}a&b\\ c &d \end{array}\right)\cdot (z+t\jj) \ = \ (a(z+t\jj)+b)(c(z+t\jj)+d)^{-1}\ .$$
Alternatively,  and more explicitly, the action is given by
$$\left(\begin{array}{ll}a&b\\ c &d \end{array}\right)\cdot (z+t\jj) \ = \
\frac{(az+b)\overline{(cz+d) } + a\overline c t^2}{|cz +d|^2 + |c|^2t^2} \ +\
\frac{t}{|cz+d|^2+|c|^2t^2}\, \jj
      \ .$$

\begin{figure}[h]
$$\includegraphics[height=9cm]{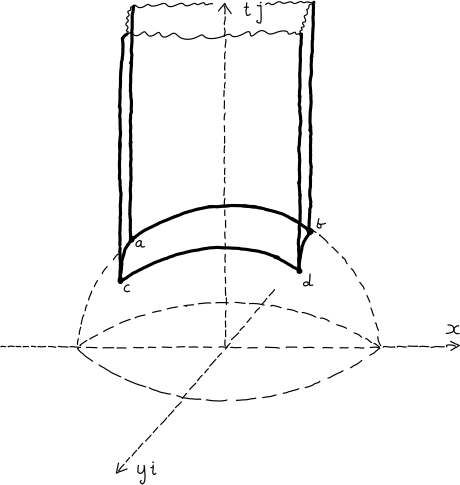}$$
\caption{Portion of the non-compact fundamental domain for the action of $SL_2({\cal O}_{-1})$ on ${\frak h}^3$}\label{FIGpicdomain}
\end{figure}

Let $G=SL_2(\mathbb Z[\ii])=SL_2({\cal O}_{-1})$.
A standard `fundamental domain' $D$ for the  restricted  action of $G$ on
${\frak h}^3$ is the  region 
\begin{equation}
D= \{z+tj\in {\frak h}^3\ |\ 0 \le |{\rm Re}(z)| \le \frac{1}{2}, 0\le {\rm Im}(z) \le \frac{1}{2}, z\overline z +t^2 \ge 1\}
\label{EQpicardregion}
\end{equation}
shown in Figure \ref{FIGpicdomain}
with some boundary points removed if one wants to minimize potential  intersections $D\cap gD$, $g\in G$ of measure zero. 
The four bottom vertices of $D$ are
$a = -\frac{1}{2}  +\frac{\sqrt{3}}{2}\jj$,
$b = \frac{1}{2}  +\frac{\sqrt{3}}{2}\jj$,
$c = -\frac{1}{2} +\frac{1}{2}\ii +\frac{\sqrt{2}}{2}\jj$,
$d = \frac{1}{2} +\frac{1}{2}\ii +\frac{\sqrt{2}}{2}\jj$. See for instance \cite[page 58]{MR1937957}.

The upper-half space ${\frak h}^3$ can be retracted onto a $2$-dimensional subspace ${\cal T} \subset {\frak h}^3$, with  ${\cal T}$
 a contractible
$2$-dimensional regular CW-complex, and where the action of  $G$ on ${\frak h}^3$
restricts to a cellular action of $G$ on ${\cal T}$. Under the restricted action there is one orbit of $2$-cells in $\cal T$, represented by the curvilinear square with vertices $a$, $b$, $c$ and $d$ in the picture. This $2$-cell has 
cyclic stabilizer group  of order $4$.
 There are three orbits of $1$-cells: the edges $ac$ and $bd$ are in the same orbit with cyclic stabilizer group of order $6$; edge $ab$ has stabilizer group
isomorphic to the quaternion group $Q4$ of order $8$; edge $cd$ has stabilizer group isomorphic to a semi-direct product $C3:C4$ of order $12$.  There are two 
   orbits of $0$-cells. Vertices $a$ and $b$ belong to the same orbit with stabilier groups isomorphic to $C3:C4$. Vertices $c$ and $d$ belong to the other orbit with stabilizer groups $SL(2,3)$.

The first $n$-terms of free $\mathbb ZH$-resolutions $R^H_\ast$
of $\mathbb Z$ for each of the finite cell-stabilizer groups $H$
can be computed using the algorithm in \cite{MR2093885}. That algorithm produces an explicit contracting homotopy on $R^H_\ast$. Using Proposition \ref{PROPwall}, these stabilizer group resolutions can be combined with $C_\ast \cal T$ to form a free
$\mathbb ZG$-resolution $R^G_\ast$ of $\mathbb Z$.

\begin{example}\label{EXpicard}
The following \hap commands use an implementation of the resolution $R^G_\ast$
for $G=SL_2(\mathbb Z[\ii])$ to compute
$$\begin{array}{lcl}
H^1(G,P_{{\cal O}_{-1}}(64)) &\cong &\mathbb Z_2^{11} \oplus \mathbb Z_4 \oplus \mathbb Z_8^2 \oplus \mathbb Z_{16} \oplus \mathbb Z_{160} \oplus \mathbb Z_{320}\ ,
\medskip\\
H^2(G,P_{{\cal O}_{-1}}(64)) &\cong &\mathbb Z^2 \oplus \mathbb Z_2^{58} \oplus \mathbb Z_{10}^2
\oplus \mathbb Z_{30} \oplus \mathbb Z_{60}^3 \oplus
\mathbb Z_{81900} \oplus \mathbb Z_{163800}\\
&&\oplus \mathbb Z_{45298780162170032823378180868002600330993000}\\
&&\oplus \mathbb Z_{90597560324340065646756361736005200661986000}\ .
\end{array}$$

\begin{verbatim}
gap> R:=ResolutionSL2QuadraticIntegers(-1,3);;
gap> G:=R!.group;;
gap> M:=HomogeneousPolynomials(G,64);;
gap> C:=HomToIntegralModule(R,M);;
gap> D:=ContractedComplex(C);;
Cohomology(D,1);
Cohomology(D,2);
gap> Cohomology(D,1);
[ 2, 2, 2, 2, 2, 2, 2, 2, 2, 2, 2, 4, 8, 8, 16, 160, 320 ]
gap> Cohomology(D,2);
[ 2, 2, 2, 2, 2, 2, 2, 2, 2, 2, 2, 2, 2, 2, 2, 2, 2, 2, 2, 2, 2, 2, 2, 2, 2, 
  2, 2, 2, 2, 2, 2, 2, 2, 2, 2, 2, 2, 2, 2, 2, 2, 2, 2, 2, 2, 2, 2, 2, 2, 2, 
  2, 2, 2, 2, 2, 2, 2, 2, 10, 10, 30, 60, 60, 60, 81900, 163800, 
  45298780162170032823378180868002600330993000, 
  90597560324340065646756361736005200661986000, 0, 0 ]
\end{verbatim}
The prime factorization of the largest torsion coefficient in the abelian invariant decomposition of 
 $H^2(G,P_{{\cal O}_{-1}}(64))$ is
$$
 2^4\times 3^2\times 5^3\times 7 \times 11\times 13\times 17\times 19\times 29\times 37\times 41\times 43\times 47\times 53\times 59\times
  61\times 197\times 103979\times 44811339594403 .
$$
\end{example}

Cuspidal cohomology for $G=SL_2(\mathbb Z[\ii])$ can be defined and implemented in a fashion directly analogous to that for $SL_2(\mathbb Z)$. The basic idea
is that the non-compact
fundamental domain $D$ of
(\ref{EQpicardregion}) is homeomorphic to $[0,1]\times[0,1]\times [0,1)$ and
can  be  compactified to $\overline{D}=[0,1]\times[0,1]\times [0,1]$ in analogy with Figure \ref{FIGfd}. Then $\overline D$ becomes the fundamental domain for a CW-complex $X^{BS}$ on which  $G$ acts (non-freely) by permuting cells. Proposition \ref{PROPwall} is then used to construct the chain complex
$C_\ast(\tilde X^{BS})$ needed to apply Definition (\ref{EQcusp}). This  is not yet implemented in \happ.

\section{Integral cohomology of $SL_2({\cal O}_d)$ and other groups}

To extend the above cohomological techniques to $G=SL_2({\cal O}_d)$ we 
require
a contractible CW-complex $\cal T$ in which $G$ acts cellularly with cell stabilizers $H$ for which we can compute a free $\mathbb ZH$-resolution $R^H_\ast$.
There are two related approaches to computing such a $\cal T$, both of which are well-documented in the literature. One approach makes use of the fact that ${\frak h}^3$ is a metric space on which $G$ acts discontinuously by isometries, and focuses on constructing a Dirichlet fundamental domain
$$D(v) =\{ w\in {\frak h}^3\ :\ d(v,w) \le d(gv,w) {\rm ~for~all~} g\in G\}$$
where $v\in {\frak h}^3$ is some suitable choice of point, and $d(\,,\,)$ 
denotes the metric on ${\frak h}^3$. The domain $D(v)$ is defined as 
an intersection of half spaces $\{w\in {\frak h}^3\ :\ d(v,w)\le d(v,g\cdot v)\}$, one half space for each $g\in G$. However, 
only finitely many of the half spaces are 
actually
 needed to determine $D(v)$; 
 using Poincar\'e's theorem (see \cite{MR0297997})
 this finite intersection can be implemented on a computer and used to determine  the face lattice of $D(v)$ with any cuspidal vertices added. 
 For a detailed account of the computation of $D(v)$ the reader should consult 
  the 
  work  of Aurel Page \cite{PageMSC,MR3356030}. The reader should also consult the related work of Alexander Rahm \cite{RahmPhD,MR2769243} which is based on the notion of a {\em Ford fundamental domain} and on papers of Swan \cite{MR0284516}, Riley \cite{MR689477}, Mendoza \cite{MR611515},  Fl\"oge \cite{MR704107}. 
  A contractible $2$-dimensional regular
CW-complex $\cal T$ arises as the orbit of a
deformation retract of  the Ford domain with cuspidal vertices added;  
details of  this $2$-complex $\cal T $
for  various  groups $SL_2({\cal O}_d)$ have been computed by Rahm and
 stored as part of a library in \happ.

 A second approach to computing $\cal T$
 uses Voronoi's  theory of {\em perfect} quadratic forms. Let $S^m_{>0}$ denote the space of positive definite symmetric $m\times m$ matrices $Q$. Such a
 matrix $Q$ corresponds to an $m$-dimensional quadratic form. The {\em cone}
 $S^m_{>0}$ is  contractible and 
 $A\in GL_m(\mathbb Z)$  acts on $S^m_{>0}$ via
$$(A,Q) \mapsto AQA^t \ $$ 
where $A^t$ denotes the transposed matrix. For a matrix $Q\in S^n_{>0}$ and
column vector
$v\in \mathbb R^m$ set
$$\begin{array}{l}
Q[v]=v^tQv\ ,\medskip\\
\rho(v) = v\,v^t\in S^m_{>0}\ ,\medskip\\
{\rm min}(Q) = {\rm min}_{0\ne v\in \mathbb Z^m}Q[v]\ ,\medskip\\
{\rm Min}(Q) = \{v\in \mathbb Z^m\ :\ Q[v]={\rm min}(Q)\}\ .
\end{array}$$
A quadratic form $Q[v]$ is said to be {\em perfect} if a quadratic form 
$P[v]$ satisfies $P[v]={\rm min}(Q)$ for all $v\in {\rm Min}(Q)$ only if $P=Q$.
 For example, the quadratic form $Q[x,y]=x^2+xy+y^2$ has ${\rm min}(Q)=1$ and 
${\rm Min}(Q)=\{(1,0), (-1,0), (0,1), (0,-1), (1,-1), (-1,1)\}$ and is perfect.
\begin{theorem}[Voronoi \cite{MR1580737}] \label{THMvoronoi} There are only finitely many perfect $m$-dimensional forms $Q$ up to $GL_m(\mathbb Z)$-equivalence, and the polyhedral cells
$${\rm Dom}(Q) = \left\{\sum_{v\in {\rm Min}(Q)} \lambda_v \rho(v)\ :\ \lambda_v> 0\right\}$$
	tessellate the {\em rational closure} of $S^m_{>0}$.
\end{theorem}
\begin{theorem}[Ash \cite{MR0427490}] \label{THMash}
There is a $GL_m(\mathbb Z)$-equivariant $m \choose 2$-dimensional CW-complex ${\cal T}$ which is a deformation retract of  $S^m_{=1}$, where $S^m_{=1}$ denotes the quotient of the rational closure of $S^m_{>0}$ obtained by identifying scalar multiples.
\end{theorem}

\begin{figure}
$$\includegraphics[height=7cm]{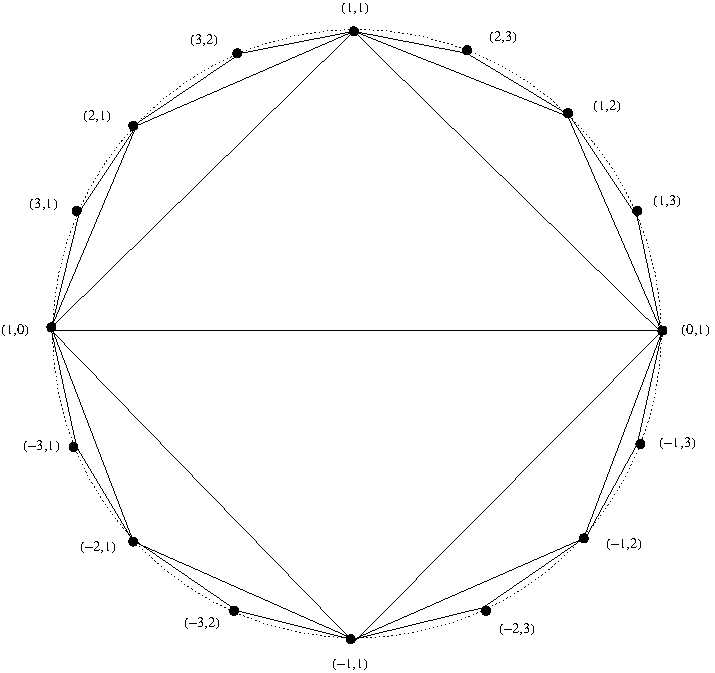}
\includegraphics[height=7cm]{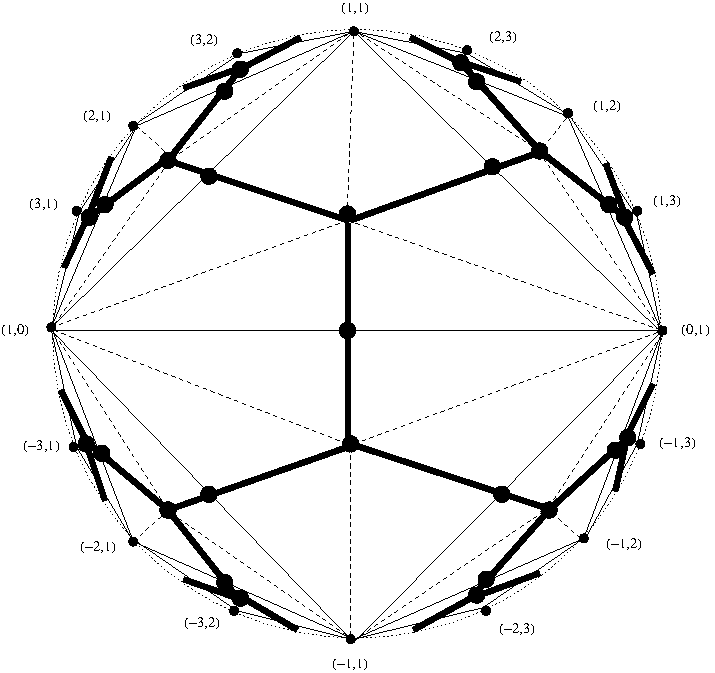}
$$
\caption{Partial Voronoi tessellation of $S^2_{=1}$ (left) and its barycentric subdivision (right). }\label{FIGtess1}
\end{figure}

It is easy to illustrate these theorems pictorially for $m=2$. The cone $S^2_{>0}$ is $3$-dimensional, and the quotient $S^2_{=1}$ of its rational closure is a $2$-dimensional regular CW-complex. We can view $S^2_{=1}$ as the union of an open unit $2$-disk with countably infinitely many points on the boundary of the $2$-disk. The cellular structure of
$S^2_{=1}$ is that of a tessellation by triangles, with the triangle vertices being the points on the boundary of the disk. A triangle vertex represents
 a ray of
quadratic forms $\mathbb R_+\rho(v)=\{\lambda \rho(v)\,:\,\lambda>0\}$ with $v\in {\rm Min}(Q)$ for some perfect form $Q$. The Voronoi tessellation
of $S^2_{=1}$ is partially pictured in Figure \ref{FIGtess1} (left), with rays $\mathbb R_+\rho(v)$ labelled simply by $v$. The barycentric subdivision of this Voronoi tessellation is partially pictured in Figure \ref{FIGtess1} (right). Some vertices and edges of the barycentric subdivision are displayed in bold. These bold vertices and edges belong to the  deformation retract of Theorem \ref{THMash} which, in the case $m=2$, is a subdivision of the cubic tree (each edge of the cubic tree is subdivided into two edges).

\begin{example}
The following \hap commands use a $3$-dimensional CW-complex $\cal T$ furnished by the theorems of Voronoi and Ash, together with
Proposition \ref{PROPwall},
 to construct a free $\mathbb ZG$-resolution $R^G_\ast$, for
 $G=SL_3(\mathbb Z)$, in degrees $\le 5$. Since
 functions specifically for
 congruence subgroups of $SL_m(\mathbb Z)$ have not yet been
implemented in \hap for $m>2$, the group $G$ is represented as a finitely presented group so that \gapp's
functionality  for finitely presented groups can be invoked. The
commands use \gapp's implementation of the low-index subgroup procedure to
list representatives of all conjugacy classes of subgroups $\Gamma \le G$ of
index at most $50$. Precisely one of these $\Gamma$ has index $48$. For this subgroup of index $48$ the commands compute
$$H_n(\Gamma,\mathbb Z) = \left\{\begin{array}{ll}
\mathbb Z_{14}, &n=1,\\
\mathbb Z_2, &n=2,\\
\mathbb Z\oplus \mathbb Z, &n=3,\\
\mathbb Z_2\oplus \mathbb Z_2\oplus \mathbb Z_2\oplus \mathbb Z_2, &n=4 .
\end{array}\right.$$
\begin{verbatim}
gap> C:=ContractibleGcomplex("SL(3,Z)");;
gap> R:=FreeGResolution(C,5);;
gap> ResolutionToResolutionOfFpGroup(R);;
gap> G:=R!.group;;
gap> index:=50;; L:=LowIndexSubgroupsFpGroup(G,index);;
gap> Length(L);
30
gap> gamma:=L[30];; Index(G,gamma);
48
gap> S:=ResolutionSubgroup(R,gamma);;
gap> C:=TensorWithIntegers(S);;
gap> D:=ContractedComplex(C);;
gap> Homology(D,1);
[ 14 ]
gap> Homology(D,2);
[ 2 ]
gap> Homology(D,3);
[ 0, 0 ]
gap> Homology(D,4);
[ 2, 2, 2, 2, 2 ]
\end{verbatim}
No contracting homotopy is implemented on the resolution $R^G_\ast$ and so
this resolution  can not yet be used to compute Hecke operators on the cohomology of $\Gamma$.
 The missing component is a contracting homotopy $h_\ast\colon C_\ast \cal T \simeq \mathbb Z$.
\end{example}

Theorems \ref{THMvoronoi}, \ref{THMash} can be extended to the case where $\mathbb Z$ is replaced by ${\cal O}_d$ and implemented on a computer as a method for determining the contractible CW-complex $\cal T$. Good accounts of this
 approach can be found, for instance, in \cite{MR2537111,MR3381447,MR2721434,MR3343219,MR3723461}. Using this approach,
Sebastian Sch\"onnenbeck has  computed a $2$-complex  $\cal T$ for  various  groups $SL_2({\cal O}_d)$ and stored its details as part of a library in \happ.
Mathieu Dutour Sikr\'ic \cite{dutour} has also used the approach to compute, and store in \happ,
higher-dimensional complexes $\cal T$ for arithmetic groups such as $SL_3(\mathbb Z[\ii])$, $SL_4(\mathbb \mathbb Z)$, $Sp_4(\mathbb Z)$.

\begin{example}
For a range of square-free values of $d$ the \hap command

\begin{verbatim}
R:=ResolutionSL2QuadraticIntegers(d,n);;
\end{verbatim}

\noindent returns $n$ dimensions of a free $\mathbb ZG$-resolution $R^G_\ast$
for $G=SL_2({\cal O}_d)$. The \hap session of Example \ref{EXpicard} can be repeated with $d=-2$ to establish:
$$\begin{array}{lcl}
H^1(SL_2({\cal O}_{-2}),P_{{\cal O}_{-2}}(64))
&\cong
&\mathbb Z^2 \oplus \mathbb Z_2^5 \oplus \mathbb Z_6^6 \oplus \mathbb Z_{12} \oplus \mathbb Z_{24}^2 \oplus \mathbb Z_{48} \oplus \mathbb Z_{96} \oplus \mathbb Z_{192}\ ,\medskip\\
H^2(SL_2({\cal O}_{-2}),P_{{\cal O}_{-2}}(64)) &\cong \end{array}$$
$$\begin{array}{l}
\mathbb Z^2 \oplus \mathbb Z_2^{34} \oplus \mathbb Z_{6}^3
\oplus \mathbb Z_{12}^9 \oplus \mathbb Z_{36}^3 \oplus
\mathbb Z_{72}^6 \oplus \mathbb Z_{144}
\oplus \mathbb Z_{4752}^4
\oplus \mathbb Z_{3792096}^2 \oplus \mathbb Z_{9347516640} \oplus \mathbb Z_{18695033280}\\
\oplus \mathbb Z_{8223545796645304770924605527348196650673670016543148734143443796390563963212512724652966920933440}^2\ .\end{array}
$$
The prime factorization of the largest torsion coefficient in the abelian invariant decomposition of
 $H^2(SL_2({\cal O}_{-2}),P_{{\cal O}_{-2}}(64))$ is
$$ 2^6 \times  3^5 \times  5\times  7\times 11\times 17\times 19\times 23\times 29\times 31\times 37\times 41\times 43\times
  47\times 53\times 59\times 61\times 138493 $$ $$\times 1367917218822877368259426449806293
\times
  1856200299217477154598445936975567 \ .$$
\end{example}

\section{Integral cohomology of $PSL_2({\cal O}_d)$ and $GL_2({\cal O}_d)$}

The contractible CW-complex $\cal T$ used in the construction of a free 
resolution for $SL_2({\cal O}_d)$ can also be used, in a similar fashion, 
to construct a resolution for $PSL_2({\cal O}_d)$. One just needs to quotient each of the
 finite slabilizer groups in $SL_2(2,{\cal O}_d)$ by the group $\langle -I\rangle$ with $I$ the identity matrix.
For instance, the  \hap commands

\begin{verbatim}
gap> R:=ResolutionPSL2QuadraticIntegers(-11,3);;
gap> M:=HomogeneousPolynomials(R!.group,5,5);;
gap> C:=HomToIntegralModule(R,M);;
gap> Cohomology(C,2);
[ 2, 2, 2, 2, 2, 2, 2, 2, 60, 660, 660, 660, 0, 0, 0, 0, 0, 0 ]
\end{verbatim}

\noindent establish

\begin{equation}
H^2(PSL_2({\cal O}_{-11}),P_{{\cal O}_{-11}}(5,5)) =
(\mathbb Z_2)^8 \oplus \mathbb Z_{60} \oplus (\mathbb Z_{660})^3 \oplus \mathbb Z^6     
\label{EQsengun}
\end{equation}
with coefficient module
$$ P_{{\cal O}_{-d}}(k,\ell) = P_{{\cal O}_{-d}}(k) \otimes_{{\cal O}_{-d}} \overline{P_{{\cal O}_{-d}}(\ell)}  $$ 
where the bar denotes a twist in the action obtained from complex conjugation. For an action of the projective linear group we must insist that $k+\ell$ is even.
The calculation (\ref{EQsengun})   was first made by Mehmet Haluk Sengun 
in \cite{MR2859903} where he records many cohomology computations for Euclidean Bianchi groups $PSL_2({\cal O}_d)$ , $d=-1,-2,-3,-7,-11$. 

For an example involving a non-Euclidean Bianchi group, the above commands can be varied to calculate
$$H^2(PSL_2({\cal O}_{-6}),P_{{\cal O}_{-6}}(64)) \cong \mathbb Z^4 \oplus A$$
where $A$ is a finite abelian group of order equal to a $1429$-digit integer; the invariant factor decomposition of  $A$ is a direct sum of 
$158$ finite cyclic groups, the largest cyclic group having order equal to a $558$-digit integer. \gapp's standard integer factorization routines 
are unable  to determine the prime decomposition of this $558$-digit integer in reasonable time.

A free resolution for $GL_2({\cal O}_d)$ can be constructed  using the short exact sequence
\begin{equation}
SL_2({\cal O}_d) \rightarrowtail GL_2({\cal O}_d) \stackrel{\det}{\twoheadrightarrow} U({\cal O}_d)
\end{equation}
in which $U({\cal O}_d)$ denotes the group of units of ${\cal O}_d$. 
When $d$ is square-free negative the group
 $U({\cal O}_d)$ is finite  of order $4$ if $d=-1$, order $6$ id $d=-3$, and order $2$ otherwise. When $d$ is square-free positive the group
 $U({\cal O}_d)$ is isomorphic to $C_2\times C_\infty$.
 We can thus construct free resolutions $R^{SL_2({\cal O}_d)}_\ast$
and $R^{U({\cal O}_d)}_\ast$ for the kernel and image of the determinant homomorphism.
 The group $GL_2({\cal O}_d)$ acts on the chain complex $R^{U({\cal O}_d)}_\ast$ in a way that each element of $R^{U({\cal O}_d)}_\ast$ is stabilized by $SL_2({\cal O}_d)$. Proposition \ref{PROPwall} can be used to construct the required free resolution $R^{GL_2({\cal O}_d)}_\ast$. This is implemented in \happ.

\section{Congruence subgroups of $SL_2({\cal O}_d)$}\label{SECod}
For a square-free integer $d$ the field $K_d=\mathbb Q(\sqrt{d})$ can be constructed as a vector space of dimension $2$ over $\mathbb Q$ endowed with a multiplication. This construction, together with conjugation
$$\overline{} K_d\colon \rightarrow K_d, a+b\sqrt{d}\mapsto a-b\sqrt{d}\ ,$$ the trace function
$${\rm tr}\colon K_d\rightarrow \mathbb Q, \alpha\mapsto \alpha  +\overline{\alpha} $$
and norm
$${\rm N}\colon K_d^\times \rightarrow \mathbb Q^\times, \alpha \mapsto \alpha \overline{\alpha}$$
are readily implemented on a computer.
 An element of $K$ is an {\em integer} if its minimal monic polynomial over 
$\mathbb Q$ has coefficients in $\mathbb Z$.  
The ring of integers ${\cal O}_{d}$ is readily implemented as a free abelian subgroup ${\cal O}_{d}=\mathbb Z \oplus \omega \mathbb Z \subset K_d$
endowed with the same multiplication, where 
$$\omega = \left\{\begin{array}{lll}
\sqrt{d}  & {\rm if~} d \equiv 2,3 {\rm ~mod~} 4\, ,\\
\frac{1+\sqrt{d}}{2}
&{\rm if~} d \equiv 1 {\rm ~mod~} 4\ .\end{array}\right.$$

An ideal $\frak{a} \triangleleft {\cal O}_d$ can be specified by giving any finite  set that generates it as an ideal. A Hermite Normal Form  
algorithm can be used to test whether an element $\alpha \in {\cal O}_d$ belongs to $\frak{a}$. It can also be used to implement addition and multiplication in the (finite) quotient ring ${\cal O}_d/\frak{a}$. 
The norm of an ideal
$$N(\frak{a}) = |{\cal O}_d/\frak{a}|$$
can be defined as the number of elements in the  quotient ring
${\cal O}_d/\frak{a}$ and can be determined from ideal generators again
using a Hermite Normal Form algorithm.
Since $N({\frak a}{\frak a}') = N({\frak a})N({\frak a}')$ an ideal is prime if its norm is a prime number. Conversely, an ideal is prime only if its norm is a prime $p$ or prime square $p^2$. Standard theory involving the quadratic character $\chi_{K_d}$ can be used to determine whether an ideal of norm $p^2$ is prime.

For any ideal $\frak a\triangleleft {\cal O}_d$ there is a canonical
  homomorphism $\pi_{\frak a}\colon SL_2({\cal O}_d) \rightarrow 
SL_2({\cal O}_d/{\frak a})$. A subgroup $\Gamma \le SL_2({\cal O}_d)$ is said 
to be a {\em congruence subgroup} if it contains $\ker\,\pi_{\frak a}$. 

The $2$-complex ${\cal T}_d$ can be used to determine generators for $G=SL_2({\cal O}_d)$ in a  fashion similar to how the cubic tree ${\cal T}_0$ was used to determine generators for $SL_2(\mathbb Z)$. Let $\Cay(G)$ denote the
 Cayley graph of $G$ with respect to these generators. The action of $G$ on $\Cay(G)$ restricts to an action of a congruence subgroup $\Gamma\le G$ of level $\frak a$  on  $\Cay(G)$. The ideal membership test for $\frak a$ can be used to implement a membership test for the group $\Gamma$, and this in turn can be used to compute a fundamental domain for the action of $\Gamma$ on $\Cay(G)$. The vertices of this fundamental domain correspond to the cosets of $\Gamma$ in $G$. The fundamental domain for $\Gamma$ can be used, for instance, 
to determine a generating set for $\Gamma$, the index of $\Gamma$ in $G$,
 and a permutation action of $G$ on the cosets of $\Gamma$. In the case when $\Gamma =\Gamma_0({\frak a})$ with $\frak a$ prime, a version of Proposition \ref{PROPcre} can be used to perform these tasks more efficiently.

\begin{example}\label{EXpenultimate}
The following \hap commands construct the prime
ideal ${\frak a}\triangleleft {\cal O}_{-1}$ in the Gaussian integers generated by the element $41+56\ii$, and then construct the congruence subgroup
$\Gamma_0({\frak a})$ of index $4818$.  

\begin{verbatim}
gap> K:=QuadraticNumberField(-1);
GaussianRationals
gap> OK:=RingOfIntegers(K);
O(GaussianRationals)
gap> a:=QuadraticIdeal(OK,41+56*Sqrt(-1));
ideal of norm 4817 in O(GaussianRationals)
gap> gamma:=HAP_CongruenceSubgroupGamma0(a);
<group of 2x2 matrices in characteristic 0>
gap> IndexInSL2O(gamma);
4818
\end{verbatim}

\noindent A maximal tree in the fundamental domain for the action of $\Gamma_0({\frak a})$ on $\Cay(SL_2({\cal O}_{-1}))$ is shown in Figure \ref{FIGbigtree}.
\begin{figure}[h]
\centerline{\includegraphics[height=4cm]{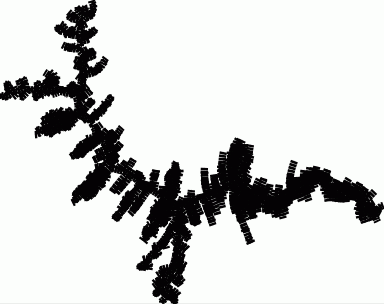}}
\caption{A maximal tree in a fundamental domain for  $\Gamma_0((41+56\ii))$. }\label{FIGbigtree}
\end{figure}
\end{example}

Our free $\mathbb ZG$-resolution $R^G_\ast$ for $G=SL_2({\cal O}_{d})$
can also be used as a free 
$\mathbb Z\Gamma$-resolution. Once the contracting homotopy of Section 
\ref{SECch} is implemented on $R^G_\ast$ ($d<0$), the resolution could be used to compute Hecke operators on the integral cohomology of $\Gamma$ using
functions currently implemented in \happ.

If one is interested only in first integral homology then the  explicit construction of a free $\mathbb Z\Gamma$-resolution can be avoided. 
 One can work  instead with a free presentation
of $\Gamma$ obtained by applying \gapp's efficient implementation of the
 Reidemeister-Schreier algorithm to a presentation of $G$.  The isomorphism
 $\Gamma^{ab}=H_1(\Gamma,\mathbb Z)$  yields the desired  homology group.
\begin{example}
The following continuation of the \hap commands of Example \ref{EXpenultimate}
establish
$$\begin{array}{lcl}
H_1(\Gamma_0({\frak a}),\mathbb Z) &\cong
&\mathbb Z_2^2\oplus\mathbb Z_4 \oplus \mathbb Z_5 \oplus \mathbb Z_7\oplus \mathbb Z_{ 16} \oplus \mathbb Z_{ 29} \oplus \mathbb Z_{43} \oplus \mathbb Z_{157}
\oplus \mathbb Z_{ 179}\\
&& \oplus \mathbb Z_{ 1877}
\oplus \mathbb Z_{ 7741} \oplus \mathbb Z_{ 22037} \oplus \mathbb Z_{292306033} 
\oplus \mathbb Z_{ 4078793513671}\end{array}$$
for the ideal ${\frak a}\triangleleft {\cal O}_{-1}$  generated by  $41+56\ii$.

\begin{verbatim}
gap> H1:=AbelianInvariants(gamma);
[ 2, 2, 4, 5, 7, 16, 29, 43, 157, 179, 1877, 7741, 22037, 292306033, 
  4078793513671 ]
\end{verbatim}
The initial terms of a free $\mathbb Z\Gamma_0({\frak a})$-resolution 
can be used to compute
$$H_2(\Gamma_0({\frak a}),\mathbb Z) \cong \mathbb Z_2\oplus \mathbb Z_2 \oplus \mathbb Z\ .$$

\textcolor{black}{For the Gaussian integers ${\cal O}_{-1}$  Conjecture (\ref{EQbergeron}) can be rewritten 
\begin{equation} \frac{\log |\Gamma_0(\aaa)_{tors}^{ab}|}{{\rm Norm}(\aaa)} \rightarrow \frac{\lambda}{18\pi},\ \lambda =L(2,\chi_{\mathbb Q(\sqrt{-1})}) = 1 -\frac{1}{9} + \frac{1}{25} - \frac{1}{49}  + \cdots \label{EQbergeron2} \end{equation} as the norm of the prime ideal $\aaa\triangleleft {\cal O}_{-1} $ tends to $\infty$. Here the value of $\lambda$ is given in terms of the $L$-function $$L(s,\chi_{\mathbb Q(\sqrt{-1})}) =\sum_{n=1}^\infty \chi_{\mathbb Q(\sqrt{-1})}(s)n^{-s}$$ of the quadratic character $\chi_{\mathbb Q(\sqrt{-1})}$ associated to the quadratic field $\mathbb Q(\sqrt{-1})$.
The equivalence between (\ref{EQbergeron}) and (\ref{EQbergeron2})  is obtained from the Humbert volume formula $$ {\rm Vol} ( {\frak h}^3 / SL_2( {\cal O}_{-d} ) ) = \frac{|D|^{3/2}}{24} \zeta_{ \mathbb Q( \sqrt{-d} ) }(2)/\zeta_{\mathbb Q}(2)                          $$ valid for square-free $d>0$, where $D$ is the discriminant of $\mathbb Q(\sqrt{-d})$, and the quadratic reciptocity formula $$\zeta_{ \mathbb Q( \sqrt{-d} ) }(s)= L(s,\chi_{\mathbb Q(\sqrt{-d})}) \zeta_{\mathbb Q}(s) $$ expressing the Dedekind zeta function as a product of an $L$-function and the Riemann zeta function. } 
The following  commands approximate the quantities
 $\lambda/18\pi = 0.0161957$ and
$\frac{\log |\Gamma_0({\frak a})_{tors}^{ab}|}{{\rm Norm}({\frak a})} = 0.0210325$ in (\ref{EQbergeron2}).

\begin{verbatim}
gap> Lfunction(K,2)/(18*3.142);
0.0161957
gap> Log_e10:=0.434294481903;; #Log_10(e)
gap> 1.0*Log(Product(H1),10)/(Norm(a)*Log_e10);
0.0210325
\end{verbatim}
\end{example}

\section{Contracting homotopies}\label{SECch}
The {\sc HAP} package \cite{hap} contains functions for computing the  Hecke operators 
$T_g\colon H^n(\Gamma,A) \rightarrow H^n(\Gamma,A)$ arising from
 any group $G$ with finite index subgroup $\Gamma <G$ and element $g\in G$ for which $\Gamma'=G\cap g\Gamma g'$ is also of finite index in $G$, and  any $\mathbb ZG$-module $A$ that is finitely generated as an abelian group. These
 functions require $n+1$ terms of a free $\mathbb ZG$-resolution
$R_\ast^G$ of $\mathbb Z$ endowed with a contracting chain homotopy.
For the groups $G=SL_2({\cal O}_d), SL_m(\mathbb Z)$ ($m\le 4$) the one ingredient that still needs to be implemented in {\sc HAP}
 is a contracting homotopy on $R^G_\ast$. In this final section we describe an approach to implementing such a contracting homotopy.

It is convenient to recall the following notion.

\begin{figure}[h]
\centerline{\includegraphics[height=4cm]{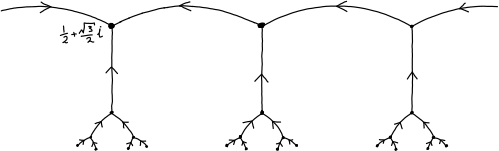}}
\caption{A portion of the cubic tree endowed with an admissible discrete vector field containing a single critical cell. Arrows $e^0_i\rightarrow e^1_{i'}$ are represented by arrow heads on the cell $e^1_{i'}$.}\label{FIGcubictreedvf}
\end{figure}

\begin{definition}\label{DEFdiscretevectorfield}
A {\em discrete vector field} on a regular CW-complex $X$ is a collection of pairs $(s , t)$, which we call {\em arrows} and denote by $s \rightarrow t$,
 satisfying
\begin{enumerate}
\item  $s,t$ are cells of $X$ with $\dim(t) =\dim(s) +1$ and with
$s$ lying in   the boundary of $t$. We say that $s$ and $t$ are {\it involved}
 in the arrow, that $s$ is the {\it source} of the arrow, and that
$t$ is the {\it target} of the arrow.

\item any cell is involved in at most one arrow.

\end{enumerate}

\end{definition}
The term {\em discrete vector field} is due to  \cite{MR1612391}.
 In an earlier work \cite{MR968920} Jones calls a very related concept  a
 {\it marking}.
By a {\it chain}\index{chain} in a discrete vector field we mean a sequence of
arrows
$$\ldots , s_1 \rightarrow t_1, s_2 \rightarrow t_2, s_3
\rightarrow t_3, \ldots $$  where the cell $s_{i+1}$
lies in the boundary of $t_i$ for each $i$.
A chain is  a {\em circuit} \index{circuit} if it is of finite length with source $s_1$
of the  initial arrow $s_1\rightarrow t_1$ lying in the boundary of  the
target $t_n$ of the final arrow $s_n\rightarrow t_n$.
A discrete vector field is
said to be {\it admissible} \index{admissible} if it contains no circuits and
  no chains that extend
infinitely to the right.
 We say that an
admissible discrete vector field is {\it maximal}\index{maximal discrete vector field} if it is not possible to
add an arrow while retaining admissibility. A cell in $X$ is said to be
{\it critical} \index{critical cell} if it is not involved in any arrow.
See Figure \ref{FIGcubictreedvf} for an example of a maximal
discrete vector field on the cubic tree, involving just one critical cell.

\begin{theorem}\cite{MR1939695,MR1612391}\label{THMbasic}
 If $X$ is a regular CW-complex with admissible discrete vector field then there is a homotopy equivalence
$$X\simeq Y$$
where $Y$ is a (possibly non-regular) CW-complex whose cells are in one-one correspondence with
  the critical cells  of $X$.
\end{theorem}

An arrow on $X$ can be viewed as representing a simple homotopy collapse, as introduced in \cite{MR0035437}. The theorem just says that an admissible discrete vector field represents some sequence of simple homotopy collapses statring at $X$ and ending at $Y$.
 At the level of cellular
chain complexes, an admissible discrete vector field on $X$
induces homomorphisms $h_{n-1}\colon C_{n-1}X \rightarrow C_nX$ of free abelian groups, defined recursively on
free generators by
$$h_{n-1}(e^{n-1}_i) = \left\{\begin{array}{ll}
0\, , &{\rm if~} e^{n-1}_i {\rm ~is~not~the~source~of~any~arrow},\\
e^n_{i'} +h_{n-1}(\partial_n(e^n_{i'})-e^{n-1}_i)\, , &{\rm if~}e^{n-1}_i\rightarrow e^n_{i'}
{\rm ~is~an~arrow~of~the~vector~field}.
\end{array}\right.$$

 In the particular case where $X$ has a single critical $0$-cell and all other cells of $X$ are involved in an arrow,  the homomorphisms $h_{n-1}$ constitute a contracting chain homotopy $H_\ast\colon C_\ast X \simeq \mathbb Z$.
The discrete vector field on the cubic tree pictured in Figure \ref{FIGcubictreedvf} corresponds to the contracting homotopy given in (\ref{EQcontracting}) and (\ref{EQfactor}).

\bigskip
As explained above,  Theorems \ref{THMvoronoi} and \ref{THMash}
provide an approach to constructing a contractible $m\choose 2$-dimensional
contractible CW-complex $\cal T$ on which $G=SL_m(\mathbb Z)$ acts with finite stabilizers, and from which one can attempt to calculate the cohomology of $G$.
A mathematically inelegant, but perhaps not totally impractical, approach to working with
 a contracting homotopy  $h_\ast\colon C_\ast \cal T \simeq \mathbb Z$ is to note that in any given computation  the values
$h(e^k_i)$ are needed on only finitely many free generators $e^k_i$ of $C_k\cal T$. So we could construct a {\em suitably large} tree $W^1$ in the
$1$-skeleton of $\cal T$ and consider the finite CW-subcomplex
$W\subset \cal T$ consisting of all cells in $\cal T$ whose closure conatins a vertex in the tree $W^1$. It {\em may} be that $W$ is contractible, and it {\em may} also
happen
 that \happ's algorithm for constructing a maximal discrete vector field on a
finite regular CW-complex would yield a discrete vector field on $W$
involving precisely one critical cell. When these two hypotheses are met we obtain a contracting homotopy on the finite subcomplex $C_\ast W \subset C_\ast{\cal T}$; if $W$ is large enough then the vector field would suffice for the computation of Hecke operators.

\begin{example}
The $3$-dimensional $G$-equivariant space $\cal T$ for $G=SL_3(\mathbb Z)$ has
one orbit of $k$-cells for $k=0,1,3$ and two orbits of $2$-cells. A full description of $\cal T$ can be found in \cite{MR0470141}.
Starting at the identity vertex $e^0\in \cal T$ of this particular $\cal T$
 and applying $n=15$ iterations of a breadth-first search, the author constructed a tree $W^1$ with $15548$ vertices;
the corresponding $3$-dimensional  CW-complex $W$ had a total of $72267$ cells. \happ's algorithm for constructing
 maximal discrete vector fields produced one on  $W$ for which there was a
single critical cell.
This contracting discrete vector field on $W$ can be viewed as a discrete
vector field on $\cal T$ which contracts the subspace $W$. The author has not yet tried to compute Hecke operators using such a discrete vector field.
\end{example}

\begin{figure}
$$\includegraphics[height=7cm, angle=270]{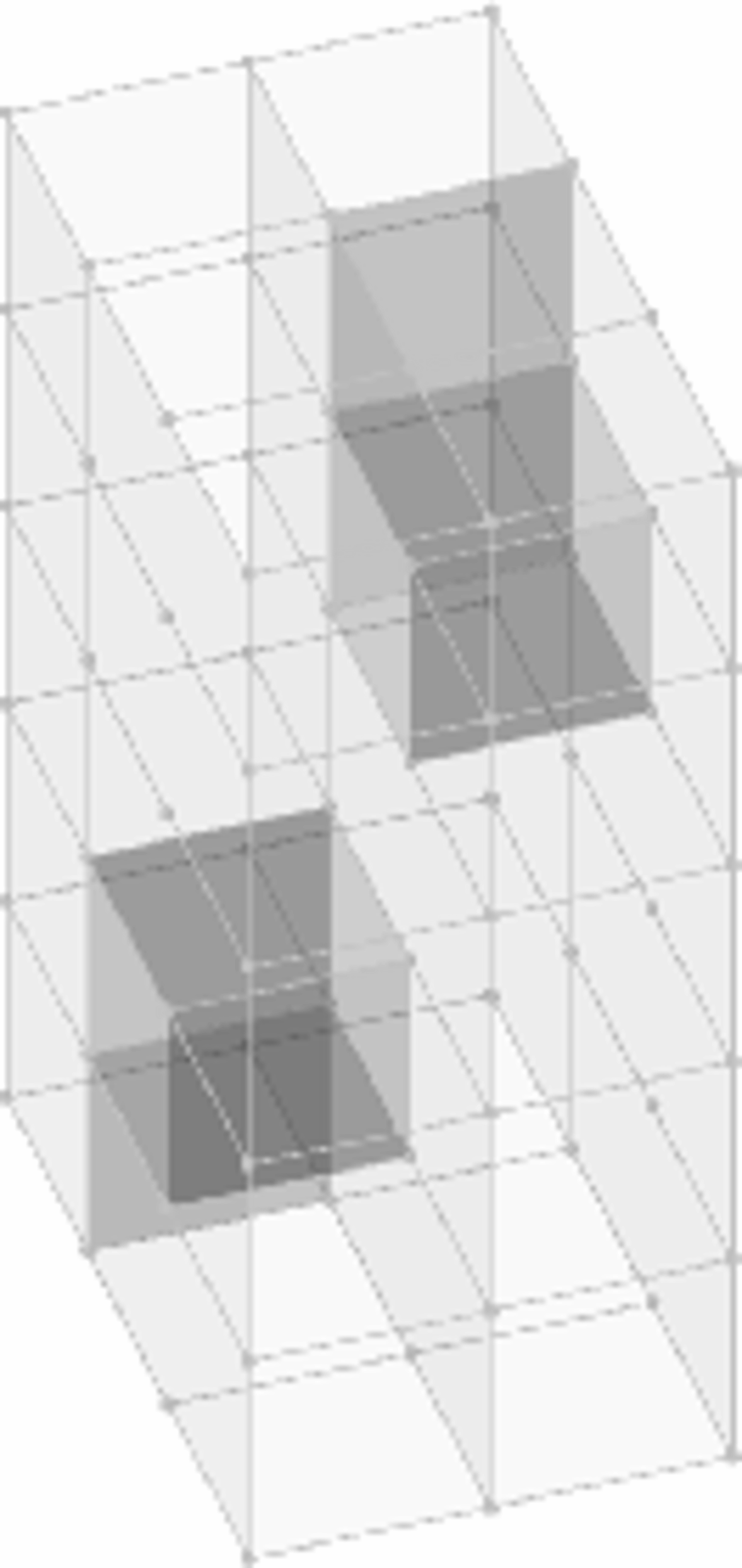}$$
\caption{Bing's house }\label{FIGbing}
\end{figure}

There are alternative approaches to constructing discrete vector fields on contractible complexes $\cal T$. Before discussing one of these, it is worth recalling that there exist contractible regular  CW-complexes that do not admit any  admissible contracting
 discrete vector field.
  Figure \ref{FIGbing} shows  a famous example, {\em Bing's house}, arising as the union of finitely many closed unit squares in $\mathbb R^3$.
The house is a
$2$-dimensional CW-complex $Y$ involving two rooms, each room having a single
entrance. The downstairs room is entered through an entrance on the roof of
the house; the upstairs room is entered through an entrance on the bottom floor
 of
the house. Suppose that Bing's house $Y$ admitted an admissible discrete vector field with precisely one critical cell $e^0$.
The arrows $e^0_i\rightarrow e^1_j$ would constitute a maximal tree in the $1$-skeleton of $Y$ rooted at the vertex $e^0$.
The remaining arrows $e^1_i\rightarrow e^2_j$ would pair  those edges not in the maximal tree with the $2$-cells of $Y$.
 Every edge $e^1_i$ is in the boundary of at least two $2$-cells, say $e^2_j$ and $e^2_{j'}$. Thus each edge $e^1_i$ which is not in the maximal tree must be   a non-initial edge in some chain $\cdots, e^k_1\rightarrow e^2_{j'}, e^1_i\rightarrow e^2_j, \cdots$ in the discrete vector field. Since the discrete vector field has only finitely many arrows it must contain a circuit. This contradicts the admissibiliy hypothesis.
To establish that Bing's house $Y$ is contractible one could use the following \hap commands to load $Y$ as a regular CW-complex involving
 72 $0$-cells, 154 $1$-cells, 83 $2$-cells,  and then compute that it is acyclic with trivial fundamental group.

\begin{verbatim}
gap> dir:=Filename(DirectoriesPackageLibrary("HAP","tst/testall")[1],"bing.txt");;
gap> Read(dir);
gap> Y:=BingsHouse;
Regular CW-complex of dimension 2
gap> Y!.nrCells(0);
72
gap> Y!.nrCells(1);
154
gap> Y!.nrCells(2);
83
gap> F:=FundamentalGroup(Y);
<fp group on the generators [  ]>
gap> Homology(Y,0);
[ 0 ]
gap> Homology(Y,1);
[  ]
gap> Homology(Y,2);
[  ]
\end{verbatim}

We now discuss an alternative approach to constructing discrete vector fields on contractible complexes.
Suppose that $X$ is a regular $n$-dimensional CW-complex, $n\ge 1$, satisfying each of the following hypotheses:
\begin{enumerate}
\item
 $X$ is {\em pure}, by which we mean
that every cell of dimension $<n$ lies in the closure of at least one $n$-cell.
\item $X$ is a  subspace of some Euclidean space $\mathbb R^n$,
	with the closure of every $n$-cell of $X$  a convex polytope (with face lattice equal to that of the polytope).
\item  $X$ is  {\em star-like},
by which we mean that there is some preferred $0$-cell $e^0\in X$ such that for any point $x\in X$ the line from $x$ to $e^0$ lies entirely in $X$.
  \end{enumerate}
For example, the  CW-complex $S^m_{n=1}$
of Theorem \ref{THMash} can be viewed as  a regular CW-complex $X$ satisfying these hypotheses.

\begin{figure}[h]
\centerline{\includegraphics[height=8cm]{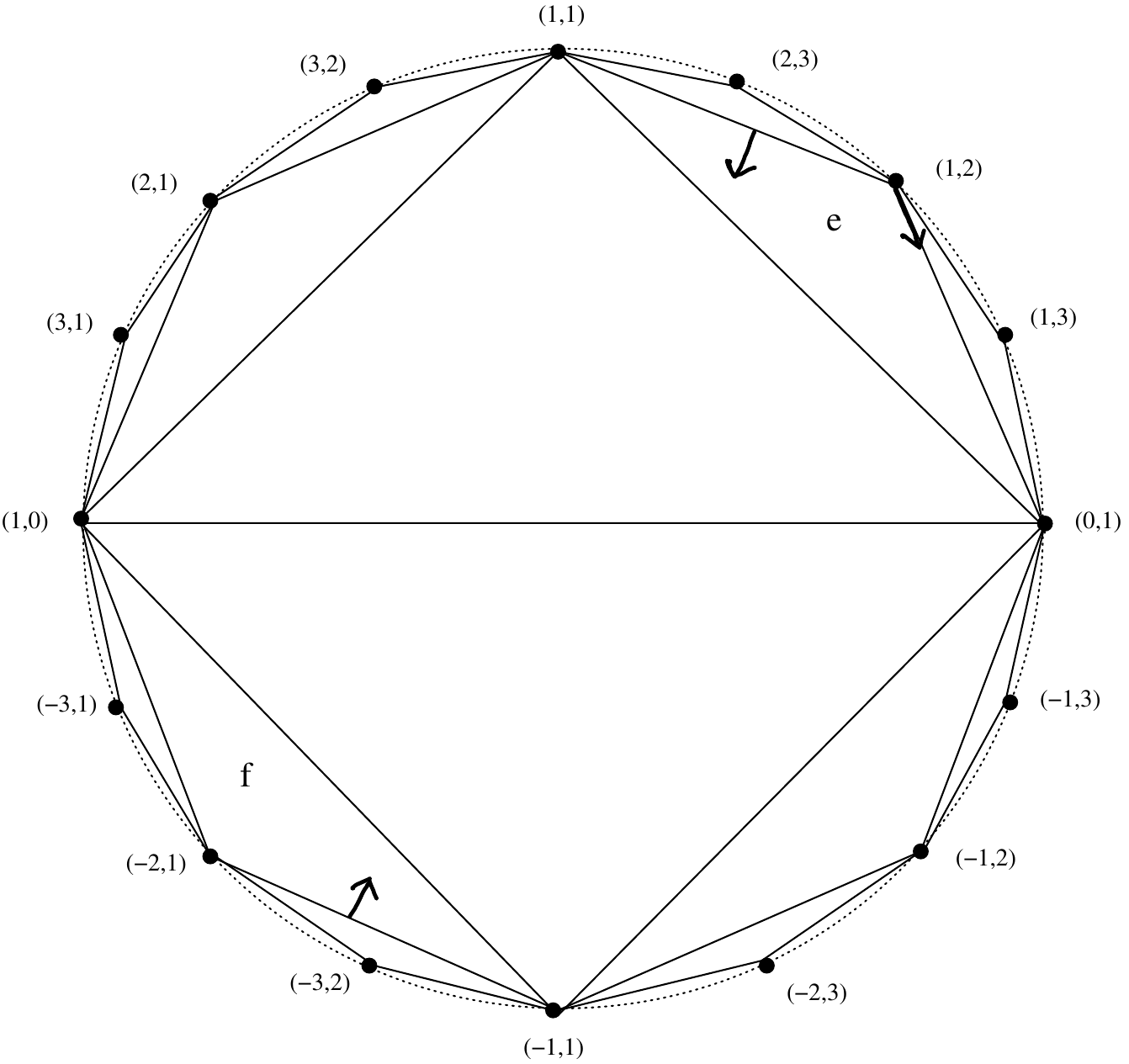}}
\caption{A portion of a contracting discrete vector field on $S^2_{=1}$}\label{FIGvoronoitessdvf}
\end{figure}

Let $X$ be  any space satisfying the  hypotheses 1--3.
For $x\in X$ let
$[x,e^0]$ denote the closed line segment from $x$ to the preferred $0$-cell $e^0$.
 We denote the closure of a $k$-cell $e^k$  by $\overline {e^k}$.
We define the {\em shadow} of a $k$-cell $e^k$ to be the set
$${\rm Sh}(e^k) =\{ x\in \overline {e^k}\ :\ [x,e^0]\cap \overline{e^k} =\{x\}~{\rm or}~ [x,e^0]\cap e^k=\emptyset \}\ .$$
The shadow ${\rm Sh}(e^k)$ is  a sub CW-complex of the closure $\overline {e^k}$ and moreover a deformation retract of $\overline {e^k}$.
Let us suppose that for each cell $e^k$  an admissible discrete vector field can be constructed on $\overline {e^k}$ for which the critical cells are precisely the cells in the shadow ${\rm Sh}(e^k)$.
The union of the discrete vector fields on the closures $\overline {e^k}$
then constitute a contracting discrete vector field on the space
$X=\bigcup \overline {e^k}$.
As an
illustration, Figure \ref{FIGvoronoitessdvf}
shows part of an admissible contracting discrete vector field on $X=S^2_{=1}$.
 In the figure two $2$-cells of $S^2_{=1}$ are labelled as $e$ and $f$; the preferred $0$-cell is labelled $(1,0)$,
and with respect to this choice the  shadow ${\rm Sh}(e)$  consists of one edge and two vertices, whereas the shadow ${\rm Sh}(f)$  consists of two edges and 
the single preferred $0$-cell; the figure shows discrete vector fields on the closures of $e$ and $f$; these vector fields are restrictions of a contracting discrete vector field on $S^2_{=1}$.

Each homotopy equivalence
${\rm Sh}(e^k)\hookrightarrow \overline {e^k}$ can indeed be realized as an admissible discrete vector field on $\overline {e^k}$ since the boundary of any convex polytope is shellable and the shelling can be constructed so that the facets in the shadow come first \cite[Theorem 8.12]{opac-b1085685}. A shelling of a contractible space gives rise to a discrete vector field on (or collapse of) the space \cite[Lemma 17]{MR1923215}. We have thus established the following.

\begin{proposition}\label{PROPtwotwo}
Let $G=SL_n({\cal O}_{d})$ for square-free $d<0$. In principle, a  free $\mathbb ZG$-resolution $R^G_\ast$ of $\mathbb Z$ and contracting homotopy $h_\ast \colon R^G_\ast \simeq \mathbb Z$ can  be implemented on a computer, with $R^G_k$ finitely generated for all $k\ge 0$.
Arbitrary elements $w \in R^G_k$  can be expressed and their
images $d_k(w)$, $h_k(w)$ can be uniquely determined.
\hfill{$\blacksquare$}\end{proposition}

\section{Acknowledgements} 
The author is very grateful to the referee for helpful comments and suggestions.

\bibliography{mybib}{}
\bibliographystyle{plain}

\end{document}